\documentclass[a4paper, 11pt]{article}

\usepackage{authblk}
\usepackage{verbatim}
\usepackage[flushleft]{threeparttable}
\usepackage{lineno}
\usepackage{amssymb}
\usepackage{amsthm}
\usepackage{graphicx}
\usepackage{caption}
\usepackage{subcaption}
\usepackage{booktabs}
\usepackage{tabularx}
\usepackage{amsmath}
\usepackage{dsfont}
\usepackage{mathtools}
\usepackage{float}
\usepackage{multirow}
\usepackage{color}
\restylefloat{table}
\usepackage{upgreek}
\usepackage{nicefrac}
\usepackage{xargs}                      
\usepackage[pdftex,dvipsnames]{xcolor} 
\usepackage[colorinlistoftodos,prependcaption,textsize=tiny]{todonotes}
\usepackage{xcolor}
\usepackage{multicol}
\usepackage{algorithm}
\usepackage[noend]{algpseudocode}
\makeatletter
\def\BState{\State\hskip-\ALG@thistlm}
\makeatother
\usepackage{nomencl}
\makenomenclature

\theoremstyle{plain}

\theoremstyle{definition}
\newtheorem{defn}{Definition}[section]

\theoremstyle{remark}

\newcommandx{\unsure}[2][1=]{\todo[linecolor=red,backgroundcolor=red!25,bordercolor=red,#1]{#2}}
\newcommandx{\change}[2][1=]{\todo[linecolor=red,backgroundcolor=red!25,bordercolor=red,#1]{#2}}
\newcommandx{\info}[2][1=]{\todo[linecolor=OliveGreen,backgroundcolor=OliveGreen!25,bordercolor=OliveGreen,#1]{#2}}
\newcommandx{\improvement}[2][1=]{\todo[linecolor=Plum,backgroundcolor=Plum!25,bordercolor=Plum,#1]{#2}}
\newcommandx{\thiswillnotshow}[2][1=]{\todo[disable,#1]{#2}}

\begin{document}
	
	\title{Uncertainty quantification for complex systems with very high dimensional response using Grassmann manifold variations}
	
	\author[1]{D.G. Giovanis\footnote{dgiovan1@jhu.edu}}
	\author[1]{M.D. Shields\footnote{michael.shields@jhu.edu}}
	
	\affil[1]{Department of Civil Engineering,  Johns Hopkins University, USA}
	
	\maketitle
	
	\begin{abstract}
		
		This paper addresses uncertainty quantification (UQ) for problems where scalar (or low-dimensional vector) response quantities are insufficient and, instead, full-field (very high-dimensional) responses are of interest. To do so, an adaptive stochastic simulation-based methodology is introduced that refines the probability space based on Grassmann manifold variations. The proposed method has a multi-element character discretizing the probability space into simplex elements using a Delaunay triangulation. For every simplex, the high-dimensional solutions corresponding to its vertices (sample points) are projected onto the Grassmann manifold. The pairwise distances between these points are calculated using appropriately defined metrics and the elements with large total distance are sub-sampled and refined. As a result, regions of the probability space that produce significant changes in the full-field solution are accurately resolved. An added benefit is that an approximation of the solution within each element can be obtained by interpolation on the Grassmann manifold without the need to develop a mathematical surrogate model. The method is applied to study the probability of shear band formation in a bulk metallic glass using the shear transformation zone theory.
	\end{abstract}

\section{Introduction} 
\label{S:1}

\noindent
Over the last several decades, the use of high-dimensional physics-based computational models (high-fidelity models) to assess the behavior of complex physical systems has become ubiquitous. Using models of increasing fidelity for the description of the physical system helps ensure that the requisite information to understand the physics of the system is available but may come at great computational cost. Given uncertainties in the model input or the parameters of the model itself, the problem is set in a probabilistic framework in order to quantify the uncertainty in the behavior/response of the physical system. This process, referred to as uncertainty quantification (UQ) becomes tremendously expensive with the use of high-fidelity models as it requires many repeated model simulations. 


Given these computational constraints, it is essential to carefully select the points at which the solution is evaluated. Over the past 20+ years, numerous approaches have been proposed to efficiently sample a probability space in order to reduce the required number of simulations. These methods differ in their approach ranging from those that are statistical in nature (i.e.\ employing Monte Carlo simulations) to those using so-called surrogate models or reduced-order models (ROMs). Statistical methods generally refer to those methods that sample randomly according to some variance reduction strategy that usually serves to improve the space-filling properties of the samples (e.g. Latin hypercube sampling)\cite{Shields_Zhang_16}. Surrogate modeling approaches aim to develop a simpler mathematical function that is inexpensive to evaluate and accurately approximates/interpolates the solution between sample points where the exact solution is known. These methods include the popular stochastic collocation method \cite{Babuska07} that employs generalized polynomial chaos interpolating functions \cite{Xiu_Karniadakis02} and the Gaussian process regression (kriging) model \cite{Bilionis13} that fits a Gaussian stochastic process to the solution, among others. ROMs, meanwhile, aim to preserve the governing equations while drastically reducing the number of degrees of freedom. ROM approaches include reduced-basis methods \cite{Tan2008,Rozza08} and nonlinear projection methods \cite{Farhat2008}.

In the past decade, it has been demonstrated that these methods can be made increasingly efficient by sampling adaptively. That is, a learning function is defined that uses the information available from existing simulations to inform the placement of additional samples. In the context of Monte Carlo methods, the second author has proposed to use gradients of the solution estimated from a surrogate model to sample in regions of high solution variance \cite{Shields}. Several authors have proposed adaptive stochastic collocation methods in which the hierarchical surplus or some other metric is used to refine the probability space using e.g.\ multi-element methods \cite{Foo2008,SSC2010,Bhaduri} or hierarchical sparse-grids \cite{Ma2009}. Adaptive methods have been used for kriging model construction by employing probabilistic learning functions such as those used in the AK-MCS approach \cite{Echard_11}. In the context of ROMs, adaptive sampling has been conducted by posing the sampling problem as a PDE-constrained optimization problem \cite{Tan2008}. 

Developing a useful learning function for adaptive sampling poses significant challenges as it typically requires an informative scalar metric of the solution and/or its error from the surrogate/ROM. For problems in reliability analysis for example, such metrics may be obvious as the analysis naturally admits a scalar quantity of interest (QOI). In more general problems where the solution is temporally and/or spatially varying, defining such a meaningful scalar QOI is often not intuitive. In this work we propose an adaptive sampling method that uses, as a refinement criterion, a distance metric that measures the difference in the intrinsic geometric structure of the full-field (time and space varying) solution through a nonlinear projection. In particular, the method arises by combining the multi-element probability space discretization 
of the simplex stochastic collocation method \cite{SSC2010} with ROM-based interpolation on each element as proposed by Amsallem and Farhat \cite{Farhat2008} and employing the doubly-infinite Grassmann distance of \cite{YeLim2014} as an element refinement metric.

The proposed method takes advantage of the fact that the geometric relationship between high-dimensional data can be investigated by projection onto low-dimensional manifold spaces, i.e.\ trajectories between solutions computed using high-dimensional models are contained and may be interpolated \cite{Farhat2008} in low-dimensional subspaces. More specifically, we project the solution onto a Grassmann manifold and use the features of this manifold to inform the UQ procedure. A Grassmann manifold is a smooth, measurable manifold on which distances between solution ``snapshots'' can be computed and interpolation performed \cite{Farhat2008}. Exploiting these property, we can identify transition regions (regions where the system's behavior changes) in the probability space and concentrate samples of the uncertain parameters in these regions to better resolve the stochastic system response. 

Thanks to the ability to interpolate on Grassmann manifolds, the method naturally yields an approximation of the full-field solution at any point in the probability space by simply reconstructing the full-field from its interpolated basis vectors and singular values.  This eliminates the need for mathematical/statistical surrogates (built using e.g. kriging or polynomial chaos expansions) which use arbitrary basis functions and cannot be fit to very high-dimensional quantities.

The paper is organized as follows. Section 2 provides an introduction to important concepts in differential geometry and more specifically properties of the Grassmann manifold. In this section, relations between the Grassmann and the Stiefel manifolds are discussed and definitions of the tangent space and the geodesic distance on the Grassmann provided. Extensions are presented to compute distances between subspaces of different dimensions as well. The concepts presented in Section 2 are admittedly abstract but essential and an elementary background in Riemannian geometry is a prerequisite. The reader is referred to \cite{Boothby2003} for further reading. In Section 3, a Delaunay triangulation-based method for discretization and refinement of a probability space is described. This multi-element sampling scheme is coupled with the Grassmann projection and interpolation of the high-dimensional model to constitute the proposed methodology in Section 4. Unlike Section 2, Sections 3-4 are intended to provide concrete algorithmic details that can be easily implemented by users with knowledge of concepts in linear algebra (i.e.\ we abandon notation based on subspaces in favor of matrix operations where possible). In Section 5, we utilize the proposed methodology to conduct uncertainty quantification on the full (40,000-dimensional) plastic strain field of an amorphous solid undergoing simple shear using the shear transformation zone theory of plasticity. The results highlight the accuracy, efficiency, and robustness of the proposed method. 

\section{The Grassmann manifold}
\label{S:2}

\noindent
Two manifolds of special interest, that arise in differential geometry and can be cast into a matrix-based approach for use in numerical linear algebra are the Stiefel manifold and the Grassmann manifold \cite{Absil2004,Edelman1998}. Some definitions from linear algebra used in order to comprehend the structure of these manifolds are: a $p$-plane is a $p$-dimensional subspace of $\mathbb{R}^n$, a $p$-frame is an ordered orthonormal basis of a $p$-plane, represented as a $n\times p$ orthonormal matrix.  The compact\footnote{The non-compact Stiefel manifold is the set of $n \times p $ matrices that have full rank.} Stiefel manifold denoted as $\mathcal{V}(p,n)$ can be defined as the set of all $p$-frames satisfying 

\begin{equation}\label{CompactStiefel}
\mathcal{V}(p,n) = \{\textbf{X}\in \mathbb{R}^{n\times p}:\textbf{X}^\intercal\textbf{X} = \textbf{I}_p\}.
\end{equation}
where $\textbf{I}_p$ is the $p\times p$ identity matrix. The Grassmann manifold $\mathcal{G}(p,n)$ is the set of all $p$-planes in $\mathbb{R}^n$, represented as an arbitrary $n$-by-$p$ orthogonal matrix. A point $\mathcal{X}$ on $\mathcal{G}(p,n)$ can be specified as a set of all linear combinations of a vector set and identified as an equivalence class under orthogonal transformation of Stiefel manifold, i.e.
\begin{equation}\label{quotientGr}
\mathcal{G}(p,n) = \mathcal{V}(p,n)/\mathcal{O}(p)
\end{equation} 
where $\mathcal{O}(p)$ is the orthogonal group (orthogonal $p\times p$ matrices).\footnote{$\mathcal{G}(p,n)$ is a quotient manifold quotient space of orthogonal groups $\mathcal{G}(p,n) = \mathcal{O}(n)/(\mathcal{O}(p) \times \mathcal{O}(n-p))$.} \footnote{$\mathcal{G}(p,n)$ and $\mathcal{V}(p,n)$ are smooth manifolds of dimension $p\times(n-p)$ and $n\times p-p\times(p+1)/2$, respectively. From a computational viewpoint a common approach is to use a matrix whose columns span a subspace of size $p\times n$ in order to represent a point on $\mathcal{G}(p,n)$ \cite{Edelman1998,Begelfor2006}.}
In this paper, we use the Stiefel manifold representation for the Grassmann manifold, i.e  each point $\mathcal{X}$ on the Grassmann manifold is an equivalence class defined by a point on the compact Stiefel manifold, represented by an orthonormal matrix $\boldsymbol{\Psi}\in \mathbb{R}^{n\times p}$.


\subsection{Riemannian structure of $\mathcal{G}(p,n)$ - Tangent space and geodesic path}

\noindent
Consider a point $\mathcal{X}=\mbox{span}(\boldsymbol{\Psi})$ on $\mathcal{G}(p,n)$, with $\boldsymbol{\Psi} \in \mathcal{V}(p,n)$ being the representative of the equivalence class. Because $\mathcal{G}(p,n)$ is a smooth (differentiable) manifold, at each point $\mathcal{X}$ there exists a tangent space, denoted $\mathcal{T}_{\mathcal{X}}$, defined by the derivative of a trajectory $\gamma(z)$ on the manifold. The tangent space, $\mathcal{T}_{\mathcal{X}}$, is a flat inner-product space with origin at the point of tangency, composed of subspaces $\mathcal{H}$ on the plane tangential to $\mathcal{G}(p,n)$ at $\mathcal{X}$. Given the matrix representation $\boldsymbol{\Psi}$ for the point $\mathcal{X}$, the tangent space can then be represented by
 \begin{equation}\label{eq:3}
  \mathcal{T}_{\mathcal{X}}\mathcal{G}(p,n) = \{\boldsymbol{\Gamma} \in \mathbb{R}^{n\times p}: \boldsymbol{\Psi}^\intercal \boldsymbol{\Gamma} = 0 \}
 \end{equation}
The trajectory $\gamma(z)$ is a geodesic path defined as the shortest path on $\mathcal{G}(p,n)$ between two points. Indexing $\gamma(z)$ on the unit line (i.e. $z\in[0,1]$), the two points are denoted $\gamma(0)=\mathcal{X}_0$ and $\gamma(1)=\mathcal{X}_1$. The geodesic path corresponds to a second order differential equation on $\mathcal{G}(p,n)$ and is therefore uniquely defined by its initial conditions $\gamma(0)=\mathcal{X}_0\in\mathcal{G}(p,n)$ and $\dot{\gamma}(0)=\dot{\mathcal{X}}_0=\mathcal{H}_0\in\mathcal{T}_{\mathcal{X}}$ and admits the following exponential mapping from the tangent space to the manifold
\begin{equation}
\exp_{\mathcal{X}_0}(\dot{\mathcal{X}}_0)=\mathcal{X}_1
\end{equation}

Let us consider that we are given matrices $\boldsymbol{\Psi}_0$ and $\boldsymbol{\Psi}_1$ spanning points $\mathcal{X}_0$ and $\mathcal{X}_1$ respectively on $\mathcal{G}(p,n)$. We aim to find an explicit mapping from the point $\boldsymbol{\Psi}_1$ onto the tangent space defined at $\boldsymbol{\Psi}_0$, represented by the matrix $\boldsymbol{\Gamma}$. More specifically, we aim to find $\boldsymbol{\Gamma}\in\mathcal{T}_{\mathcal{X}_0}$ such that the exponential mapping $\exp_{\mathcal{X}_0}(\boldsymbol{\Gamma})=\boldsymbol{\Psi}_1$ is defined locally. The sought-after inverse defines a logarithmic mapping as $\log_{\mathcal{X}_0}(\boldsymbol{\Psi}_1)=\boldsymbol{\Gamma}$. The purpose of this mapping is that it defines a geodesic path between $\boldsymbol{\Psi}_0$ and $\boldsymbol{\Psi}_1$, which further allows distances to be computed on the manifold (as we will see in the next section) and enables interpolation between the points \cite{Farhat2008}.

To arrive at this mapping, let us equivalently denote $\boldsymbol{\Gamma}$ through its thin-singular value decomposition (SVD), i.e $\boldsymbol{\Gamma} = \textbf{U}\mathbf{S}\textbf{V}^\intercal$ and thus: 

\begin{equation}\label{eq:4}
\boldsymbol{\Psi}_1 = \exp_{\mathcal{X}}(\textbf{U}\mathbf{S}\textbf{V}^\intercal)=\boldsymbol{\Psi}_0\textbf{V}\cos(\mathbf{S})\textbf{Q}^\intercal + \textbf{U}\sin(\mathbf{S})\textbf{Q}^\intercal
\end{equation}

where $\textbf{Q}$ is an $n \times n$ orthogonal matrix. Requiring that $\boldsymbol{\Gamma}$ is in the tangent space, i.e $\boldsymbol{\Psi}_0^\intercal \boldsymbol{\Gamma} = 0$,defines the following set of equations

\begin{subequations}
\begin{equation}\label{eq:6a}
\textbf{V}\cos(\mathbf{S})\textbf{Q}^\intercal = \boldsymbol{\Psi}_0^\intercal \boldsymbol{\Psi}_1
\end{equation}
\begin{equation}\label{eq:6b}
\textbf{U}\sin(\mathbf{S})\textbf{Q}^\intercal =\boldsymbol{\Psi}_1- \boldsymbol{\Psi}_0\boldsymbol{\Psi}_0^\intercal\boldsymbol{\Psi}_1. 
\end{equation}
\end{subequations}
Multiplying \eqref{eq:6a} by the inverse of \eqref{eq:6b} yields
\begin{equation}\label{eq:7}
\textbf{U}\tan(\mathbf{S})\textbf{V}^\intercal = (\boldsymbol{\Psi}_1- \boldsymbol{\Psi}_0\boldsymbol{\Psi}_0^\intercal\boldsymbol{\Psi}_1) (\boldsymbol{\Psi}_0^\intercal \boldsymbol{\Psi}_1)^{-1}
\end{equation}
Thus, the exponential mapping can be performed by taking the SVD of the matrix $\textbf{M} = (\boldsymbol{\Psi}_1- \boldsymbol{\Psi}_0\boldsymbol{\Psi}_0^\intercal\boldsymbol{\Psi}_1) (\boldsymbol{\Psi}_0^\intercal \boldsymbol{\Psi}_1)^{-1}=\textbf{U}\mathbf{S}\textbf{V}^\intercal$. The logarithmic map follows as \cite{Begelfor2006}:
\begin{equation}\label{eq:8}
\log_{\mathcal{X}}(\boldsymbol{\Psi}_1)=\boldsymbol{\Gamma} =\textbf{U}\tan^{-1}(\mathbf{S})\textbf{V}^\intercal
\end{equation}

\begin{figure}[!ht]
	\centering\includegraphics[width=1.\columnwidth]{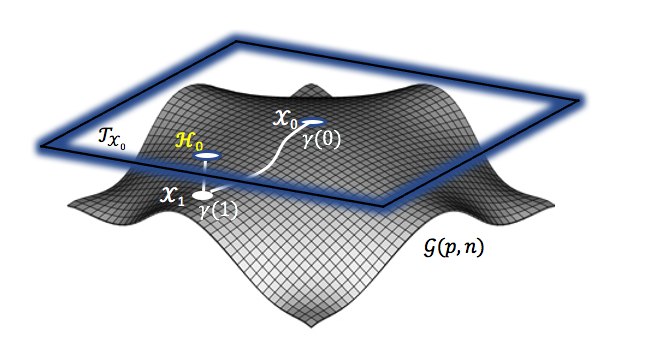}
	\caption{A Grassmann manifold $\mathcal{G}(p,n)$ with its tangent space $\mathcal{T}_{\mathcal{X}_0}$ at $\mathcal{X}_0$. A point $\mathcal{H}$ on $\mathcal{T}_{\mathcal{X}}$ can be projected to a point $\mathcal{X}_1 \in\mathcal{G}(p,n)$ through $\exp_{\mathcal{X}_0}(\mathcal{H})=\mathcal{X}_1$. The two points are connected by a geodesic path $\gamma$ of minimum length.}
	\label{Gr03}
\end{figure}

From this representation, the geodesic path $\gamma(z)$ (Fig.(\ref{Gr03})) connecting points $\mathcal{X}_0$ and $\mathcal{X}_1$ on $\mathcal{G}(p,n)$ can be expressed by \cite{Edelman1998,Boothby2003}
\begin{equation}\label{eq:9}
\gamma(z) = \mbox{span}\left[\left(\boldsymbol{\Psi}_0\textbf{V} \cos (z\mathbf{S}) + \textbf{U} \sin ( z\mathbf{S})\right)\textbf{V}^\intercal\right]
\end{equation}
In Eq. (\ref{eq:9}) the rightmost $\textbf{V}^\intercal$ can be omitted and still represent the same equivalence class as $\gamma(z)$; however, due to consistency conditions the tangent vectors used for computations must be altered in the same way and for this reason everything is multiplied by $\textbf{V}^\intercal$  \cite{Edelman1998}. 
 
\subsubsection{Metrics on the Grassmann manifold}
 
\noindent
The tangent space $\mathcal{T}_{\mathcal{X}}\mathcal{G}(p,n)$ is endowed with the canonical inner product which states that $\forall \boldsymbol{\Gamma}_1, \boldsymbol{\Gamma}_2 \in \mathcal{T}_{\mathcal{X}}\mathcal{G}(p,n): \langle \boldsymbol{\Gamma}_1,\boldsymbol{\Gamma}_2 \rangle = \mbox{tr}(\boldsymbol{\Gamma}_1^\intercal \boldsymbol{\Gamma}_2)$. This metric descends to a Riemannian metric on $\mathcal{G}(p,n)$ and in turn induces a geodesic distance. The geodesic distance between $\boldsymbol{\Psi}_0$ and $\boldsymbol{\Psi}_1$ whose columns span $\mathcal{X}_0$ and $\mathcal{X}_1 \in \mathcal{G}(p,n)$ is given by 
\begin{equation}\label{eq:10}
d_{\mathcal{G}(p,n)}(\boldsymbol{\Psi}_0,\boldsymbol{\Psi}_1) = \bigg( \sum_{i=1}^{p} \theta_i^{2} \bigg)^{\frac{1}{2}} = ||\cos^{-1}\mathbf{S}||_{F}
\end{equation}
where the principal angles $\theta_i$ can be readily computed using the full SVD of
\begin{equation}\label{eq:11}
\boldsymbol{\Psi}_0^\intercal \boldsymbol{\Psi}_1 = \textbf{U}\mathbf{S} \textbf{V}^\intercal
\end{equation}
with $\textbf{U}, \mathbf{S}, \textbf{V} \in \mathbb{R}^{p\times p}$ and the subscript $F$ denotes the Frobenius norm. The principal angles\footnote{The principal angles are bounded between $0\leq \theta_1 \leq \ldots \theta_p \leq \pi/2$} are given by $\theta_i = \cos^{-1}\sigma_i$ for $i=1,\ldots,p$ where $\sigma_i$ can be calculated from $\mathbf{S} = \mbox{diag}(\sigma_1, \ldots, \sigma_p) \in \mathbb{R}^{p \times p}$. Usually Eq.(\ref{eq:11}) is written as $ \boldsymbol{\Psi}_0^\intercal \boldsymbol{\Psi}_1 = \textbf{U}(\cos(\boldsymbol{\Theta})) \textbf{V}^\intercal$ where $\boldsymbol{\Theta} = \mbox{diag}(\theta_1, \ldots, \theta_p)$. The distance of Eq.\ (\ref{eq:10}) is called the Grassmann distance and is the primary distance used in this work. Descriptions of the various distances between subspaces considered in this work are given in Table 1 \cite{YeLim2014}.\footnote{Ye and Lim \cite{YeLim2014} show that any notion of distance between $p$-dimensional subspaces in $\mathbb{R}^n$ that depends only on the relative positions of the subspaces must be a function of their principal angles.}
\begin{table}[!ht]
\centering
	\captionof{table}{Distances on $\mathcal{G}(p,n)$ in terms of principal angles and orthonormal matrices.} 
	\begin{tabular}{c c c}   \toprule
		\textbf{Name} & $d_{\mathcal{G}(p,n)}(\boldsymbol{\Psi}_0, \boldsymbol{\Psi}_1)$ & matrix notation  \\\midrule
		Grassmann         & $\bigg( \sum_{i=1}^{p} \theta_i^{2} \bigg)^{\frac{1}{2}}$                                           & $||\cos^{-1}\boldsymbol{\Sigma}||_{F}$     \\ 
		Chordal        & $\bigg(\sum_{i=1}^k \sin^2\theta_i\bigg)^{1/2}$      & $\frac{1}{\sqrt{2}}||\boldsymbol{\Psi}_0\boldsymbol{\Psi}_0^\intercal -    \boldsymbol{\Psi}_1\boldsymbol{\Psi}_1^\intercal||_F$    \\ 
		Procrustes     & $2\bigg(\sum_{i=1}^k \sin^2(\theta_i/2)\bigg)^{1/2}$         & $||\boldsymbol{\Psi}_0\textbf{U} -    \boldsymbol{\Psi}_1\textbf{V}||_F$    \\ 
		\hline
		\label{tab:1}
	\end{tabular}
\end{table}
 
\subsection{Distance between subspaces of different dimension}
\noindent
A common problem in computing subspace distances that often arises in applications is that the subspaces of interest are of different dimension. This is the case in the proposed approach since, for different realizations of the uncertain parameters the we expect to have different behavior in the solution space. This is reflected in the rank of the matrices $\boldsymbol{\Psi}$, $\boldsymbol{\Sigma}$ and $\boldsymbol{\Phi}$ resulting from the factorization of the solution matrix $\textbf{F}$. To overcome this challenge, we rely upon an extension of the Grassmann distance between equidimensional subspaces to subspaces of different dimensions. 

To define a general distance metric between $\boldsymbol{\Psi}_i\in \mathcal{G}(r_i, n)$ and $\boldsymbol{\Psi}_j\in \mathcal{G}(r_j, n)$ where $r_i, r_j, n \in \mathbb{N}$ and $r_i \neq r_j$ one needs to isometrically embed $\mathcal{G}(r_i, n)$ and $\mathcal{G}(r_j, n)$ into an ambient Riemannian manifold and then define the distance between $\boldsymbol{\Psi}_i$ and $\boldsymbol{\Psi}_j$ as their distance in the ambient space. However, this distance depends on both the embedding and the ambient space \cite{Conway1996}. To avoid this dependence, in this work we adopt the distance proposed in \cite{YeLim2014} which is intrinsic to the Grassmaniann and independent of $n$, i.e., a $r_i$-plane $\boldsymbol{\Psi}_i$ and an $r_j$-plane $\boldsymbol{\Psi}_j$ in $\mathbb{R}^n$ will have the same distance if we consider them subspaces in $\mathbb{R}^n$ where $n \geq \min(r_i, r_j)$.  Before introducing this distance some definitions from differential geometry are required which are presented next. For their rigorous mathematical definition the reader is referred to \cite{YeLim2014}.   

\begin{defn}
The infinite Grassmannian $\mathcal{G}(r_i, \infty)$ of $r_i$-planes is defined as
\begin{eqnarray}\label{eq:12}
\mathcal{G}(r_i, \infty) = \bigcup_{n = r_i}^{\infty} \mathcal{G}(r_i, n)
\end{eqnarray}
where $\mathcal{G}(r_i, n)\subset \mathcal{G}(r_i, n+1) $ by identifying $\mathcal{G}(r_i,n)$ with $\iota_n(\mathcal{G}(r_i,n))$ where the inclusion map $\iota_n:\mathbb{R}^n\to\mathbb{R}^{n+1}, \iota_n(x_1,\dots,x_n)=(x_1,\dots,x_n,0)$. 
\label{defn:inf_Grass}
\end{defn}

Definition \ref{defn:inf_Grass} allows us to consider a matrix $\boldsymbol{\Psi}\in\mathcal{G}(p,n)$ to also exist, through isometric mapping, in all $\mathcal{G}(p,m), m>n$.

\begin{defn}
The doubly infinite Grassmannian \cite{YeLim2014} can be viewed, informally, as the disjoint union of all $r_i$-dimensional subspaces over all $r_i \in \mathbb{N}$ as 
\begin{eqnarray}\label{eq:13}
\mathcal{G}(\infty, \infty) = \coprod_{ r_i=1}^{\infty} \mathcal{G}(r_i, \infty).
\end{eqnarray}
\end{defn}

The doubly infinite Grassmannian is non-Hausdorff and therefore non-metrizable \cite{Hausdorff}. Consequently, to define a metric between any pair of subspaces of arbitrary dimensions is to define one on $\mathcal{G}(\infty,\infty)$. Many such metrics on $\mathcal{G}(\infty,\infty)$ relate poorly to the true geometry of the Grassmannian. Ye and Lim \cite{YeLim2014} have proposed a rather elegant set of distance metrics on $\mathcal{G}(\infty,\infty)$ that do, however, capture the geometry of $\mathcal{G}(k,n)\hspace{3pt}\forall k<n$ and generalize to the common distances between subspaces of equal dimension (such as and those distances in Table \ref{tab:1}). To understand this metric requires first another definition.

\begin{defn}
Let $r_i, r_j, n \in \mathbb{N}$ such that $r_i \leq r_j \leq n$. For any $\boldsymbol{\Psi}_i\in \mathcal{G}(r_i, n)$ and $\boldsymbol{\Psi}_j\in \mathcal{G}(r_j, n)$, define the subsets:
\begin{eqnarray}
\Omega_{+}(\boldsymbol{\Psi}_i) = \{\textbf{X} \in \mathcal{G}(r_j, n) : \boldsymbol{\Psi}_i \subseteq \textbf{X}\}, \quad \Omega_{-}(\boldsymbol{\Psi}_j) = \{\textbf{Y} \in \mathcal{G}(r_i, n) : \textbf{Y} \subseteq \boldsymbol{\Psi}_j\} \nonumber
\end{eqnarray}
where $\Omega_{+}(\boldsymbol{\Psi}_i)$ is called the Schubert variety \cite{Shubert} of $r_j$-planes containing $\boldsymbol{\Psi}_i$ and $\Omega_{-}(\boldsymbol{\Psi}_j)$ the Schubert variety of $r_i$-planes contained in $\boldsymbol{\Psi}_j$.
\end{defn}

The Schubert varieties $\Omega_{+}(\boldsymbol{\Psi}_i)$ and $\Omega_{-}(\boldsymbol{\Psi}_j)$ are closed subsets of $\mathcal{G}(r_j, n)$ and $\mathcal{G}(r_i, n)$ respectively that are uniquely defined by $\boldsymbol{\Psi}_i$ and $\boldsymbol{\Psi}_j$. The distance $\delta(\boldsymbol{\Psi}_i,\boldsymbol{\Psi}_j)$ measured in $\mathcal{G}(r_i, n)$ is the one between the $r_i$-plane $\boldsymbol{\Psi}_i$ and the closest $r_i$-plane $\textbf{Y}$ contained in $\boldsymbol{\Psi}_j$, i.e. is the Grassmann distance from $\boldsymbol{\Psi}_i$ to the closed subset $\Omega_{-}(\boldsymbol{\Psi}_j)$
\begin{eqnarray}\label{eq:14}
\delta(\boldsymbol{\Psi}_i,\boldsymbol{\Psi}_j)=d_{\mathcal{G}(r_i, n)}(\boldsymbol{\Psi}_i,\Omega_{-}(\boldsymbol{\Psi}_j))=\min\{d_{\mathcal{G}(r_i, n)}(\boldsymbol{\Psi}_i,\textbf{Y}):\textbf{Y}\in\Omega_{-}(\boldsymbol{\Psi}_j)\}
\end{eqnarray}

Similarly, the distance $\delta(\boldsymbol{\Psi}_i,\boldsymbol{\Psi}_j)$ measured in $\mathcal{G}(r_j, n)$ is the distance from the $r_j$-plane $\boldsymbol{\Psi}_j$ and the closest $r_j$-plane $\textbf{Y}$ containing $\boldsymbol{\Psi}_i$. Again, this is the Grassman distance $d_{\mathcal{G}(r_j, n)}$ from $\boldsymbol{\Psi}_j$ to the closed subset $\Omega_{+}(\boldsymbol{\Psi}_i)$. 
\begin{eqnarray}\label{eq:15}
\delta(\boldsymbol{\Psi}_i,\boldsymbol{\Psi}_j)=d_{\mathcal{G}(r_j, n)}(\boldsymbol{\Psi}_j,\Omega_{+}(\boldsymbol{\Psi}_i))=\min\{d_{\mathcal{G}(r_j, n)}(\boldsymbol{\Psi}_j,\textbf{X}):\textbf{X}\in\Omega_{+}(\boldsymbol{\Psi}_i)\}
\end{eqnarray}

It is shown in \cite{YeLim2014} that, for sufficiently large $n \geq r_j+1$, these two distances are equivalent 
\begin{eqnarray}\label{eq:16}
d_{Gr(r_i, n)}(\boldsymbol{\Psi}_i,\Omega_{-}(\boldsymbol{\Psi}_j)) = d_{Gr(r_j, n)}(\boldsymbol{\Psi}_j,\Omega_{+}(\boldsymbol{\Psi}_i))
\end{eqnarray}
Note that $\delta(\boldsymbol{\Psi}_i,\boldsymbol{\Psi}_j)$ is a distance from a point to set and does not define a metric on the set of all subspaces of all dimensions.  But, the equality in Eq. \eqref{eq:16} implies a distance between subspaces of different dimension with the following important properties
\begin{itemize}
\item $\delta(\boldsymbol{\Psi}_i,\boldsymbol{\Psi}_j)$ is independent of the dimension, $n$, of the ambient space and is the same for all $n\ge r_j+1$
\item $\delta(\boldsymbol{\Psi}_i,\boldsymbol{\Psi}_j)$ reduces to the Grassmann distance when $r_i=r_j$
\item $\delta(\boldsymbol{\Psi}_i,\boldsymbol{\Psi}_j)$ can be computed explicitly as
\begin{equation}\label{eq:17}
\delta(\boldsymbol{\Psi}_i,\boldsymbol{\Psi}_j) = \bigg(\sum_{k=1}^{\min(r_i, r_j)} \theta_k^2(\boldsymbol{\Psi}_i,\boldsymbol{\Psi}_j) \bigg)^{1/2}
\end{equation}
where $\theta_k(\boldsymbol{\Psi}_i,\boldsymbol{\Psi}_j)$ is the $k^{th}$ principal angle between $\boldsymbol{\Psi}_i$ and $\boldsymbol{\Psi}_j$.
\end{itemize}




To transform $\delta(\boldsymbol{\Psi}_i,\boldsymbol{\Psi}_j)$ into a metric on $\mathcal{G}(\infty, \infty)$, Ye and Lim \cite{YeLim2014} propose to set $\theta_{r_i+1} = \ldots = \theta_{r_j} = \pi/2$ (equivalent to completing $\boldsymbol{\Psi}_i$ to an $r_j$ dimensional subspace of $\mathbb{R}^n$, by adding $r_j-r_i$ vectors orthonormal to the subspace $\boldsymbol{\Psi}_j$). In this case, the distances presented in Table \ref{tab:1} can be written as shown in Table \ref{tab:2}.
\begin{table}[!ht]
\centering
	\caption{Distances on $\mathcal{G}(\infty,\infty)$ in terms of principal angles.} 
	\begin{tabular}{c c}   \toprule
		\textbf{Name} & $d_{\mathcal{G}(\infty,\infty)}(\boldsymbol{\Psi}_i,\boldsymbol{\Psi}_j)$   \\\midrule
		Grassmann         & $\bigg(|r_j-r_i|\pi^2/4 + \sum_{k=1}^{\min(r_i, r_j)} \theta_k^2 \bigg)^{1/2}$                                    \\ 
		Chordal        & $\bigg(|r_j-r_i| + \sum_{k=1}^{\min(r_i, r_j)} \sin^2\theta_k \bigg)^{1/2}$         \\ 
		Procrustes     & $\bigg(|r_j-r_i| + \sum_{k=1}^{\min(r_i, r_j)} \sin^2(\theta_k/2) \bigg)^{1/2}$       \\ 
		\hline
		\label{tab:2}
	\end{tabular}
\end{table}

The natural interpretation of $d_{\mathcal{G}(\infty, \infty)}(\boldsymbol{\Psi}_i,\boldsymbol{\Psi}_j)$ is the following \cite{YeLim2014}: the distance from $\boldsymbol{\Psi}_j$ to the furthest $r_j$-dimensional subspace that contains $\boldsymbol{\Psi}_i$, which equals the distance from $\boldsymbol{\Psi}_i$ to the furthest $r_i$-dimensional subspace contained in $\boldsymbol{\Psi}_j$. Furthermore, $d_{\mathcal{G}(\infty, \infty)}(\boldsymbol{\Psi}_i,\boldsymbol{\Psi}_j)$ can be considered as the amalgamation of two pieces of information, the distance on $\mathcal{G}(r_i,\infty)$ given by $\delta(\boldsymbol{\Psi}_i,\boldsymbol{\Psi}_j)$ and the difference in dimensions $|\mbox{dim}\textbf{X} - \mbox{dim}\textbf{Y}|$. For $r_i = r_j$ the distances of Table \ref{tab:2} reduce to the corresponding ones of Table \ref{tab:1}. 

\section{Delaunay triangulation-based mesh refinement}
\label{S:3}

\noindent
To achieve a sample design that is, in some sense, optimized for a specific operation, knowledge of the behavior of the operation over the full range of all random variables is required \cite{Shields}. Thus, an adaptive approach that seeks to identify regions of the probability space over which the solution possesses significant variability appears attractive. In this work, a multi-element technique similar to that used in the simplex stochastic collocation method (SSC) \cite{SSC2010} is used. The method utilizes a Delaunay triangulation (DT) to discretize the probability space,  and has the advantage that mesh refinement can be easily performed. Utilizing the concepts presented in Section 2, the mesh refinement is based on observed variations on the Grassmann manifold of the very high dimensional solution. The refinement criterion will be elaborated further in Section 4, while the details of how the elements are refined are discussed herein.

Consider a complete probability space $[\Omega, \mathcal{F},\mathbb{P}]$ where $\Omega$ is the sample space, $\mathcal{F}$ is the $\sigma$-algebra of the events, and $\mathbb{P}$ is a probability measure. In order to identify and locally refine regions of the multi-dimensional probability space where the operation shows variations in its behavior, an unstructured sample grid is found more attractive than the common structured tensor product or sparse grids. Let $\textbf{s}=\{s_{0},\cdots,s_{n_d-1}\} \in \Omega \subseteq \mathbb{R}^{n_d} $ be the vector of uncertain model parameters with distribution $f_\mathbf{s}$. Without loss of generality, we set $\Omega=C^{n_d}=[0,1]^{n_d}$ and denote the parameter vector by $\boldsymbol{\xi}$ with components distributed uniformly in [0,1].
By setting $\Omega=C^{n_d}$, a linear and/or nonlinear mapping is assumed to exist in order to transform from $C^{n_d}$ to any parameter domain $S$. This is strictly for convenience and is not required as the proposed method can be equivalently applied on non-hypercube probability spaces.

\begin{figure}[!ht]
	\begin{subfigure}[b]{0.5\textwidth}
		\centering\includegraphics[scale=0.23]{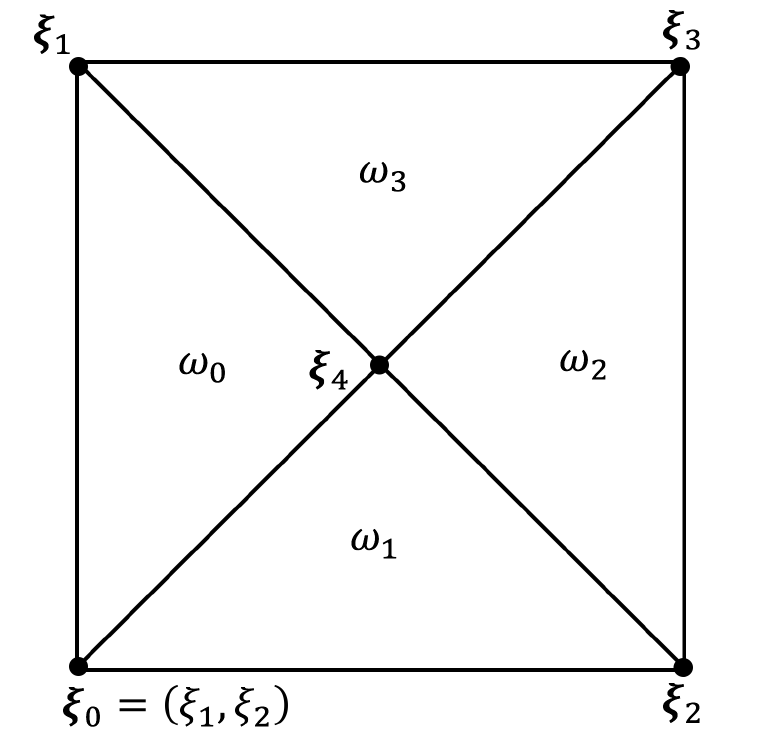}
		\caption{}
		\label{2DiniGrid}
	\end{subfigure}
	~ 
	\begin{subfigure}[b]{0.5\textwidth}
		\centering\includegraphics[scale=0.24]{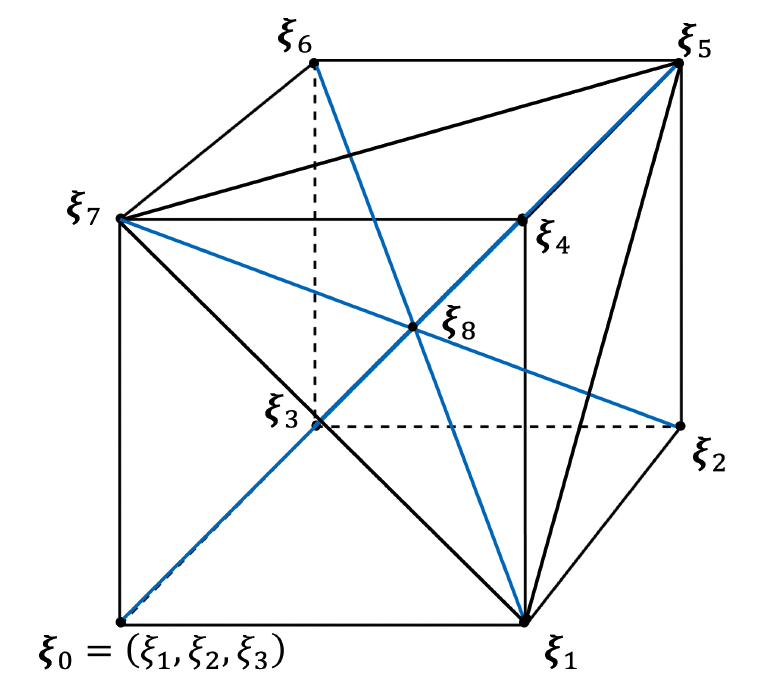}
		\caption{}
		\label{3DiniGrid}
	\end{subfigure}
	\caption{Initial sampling grid for $n_d= 2,3$: (a) $2^{n_d=2}+1= 5$ sample points and $n_{e} = 4$ triangles with $n_d+1=2+1 = 3$ vertices each, (b) $2^{n_d=3}+1 = 9$ sample points and $n_{e} = 12$ triangles with $n_d+1=3+1 = 4$ vertices each}\label{Grid}
\end{figure}

If $n_{s}$ realizations of $\boldsymbol{\xi}$ (denoted $\boldsymbol{\xi}_k,\hspace{3pt}k=0\dots n_s-1$) are available, $\Omega$ can be discretized into $n_{e}=O(n_{s}^{n_{d}/ 2})$ disjoint and space-filling elements (simplexes), i.e. $\Omega = \omega_0 \cup \ldots \cup \omega_{n_e-1}$, $\omega_i\cap\omega_j=\emptyset$, using a DT with the points $\boldsymbol{\xi}_k$ serving as the vertices.  For our purposes, the initial set consists of $2^{n_d}$ samples located at the corner points of $C^{n_d}$ and one in the interior of $C^{n_d}$ (a total of $n_{s}= 2^{n_d}+1$ points) and yielding a total of $n_e=O(n_{s}^{n_{d}/ 2})$ elements. Figs.\ (\ref{2DiniGrid}) and (\ref{3DiniGrid})  depict the initial samples and DTs for a two- and a three-dimensional problem, respectively.  

\begin{figure}[!ht]
	\centering\includegraphics[width=0.5\columnwidth]{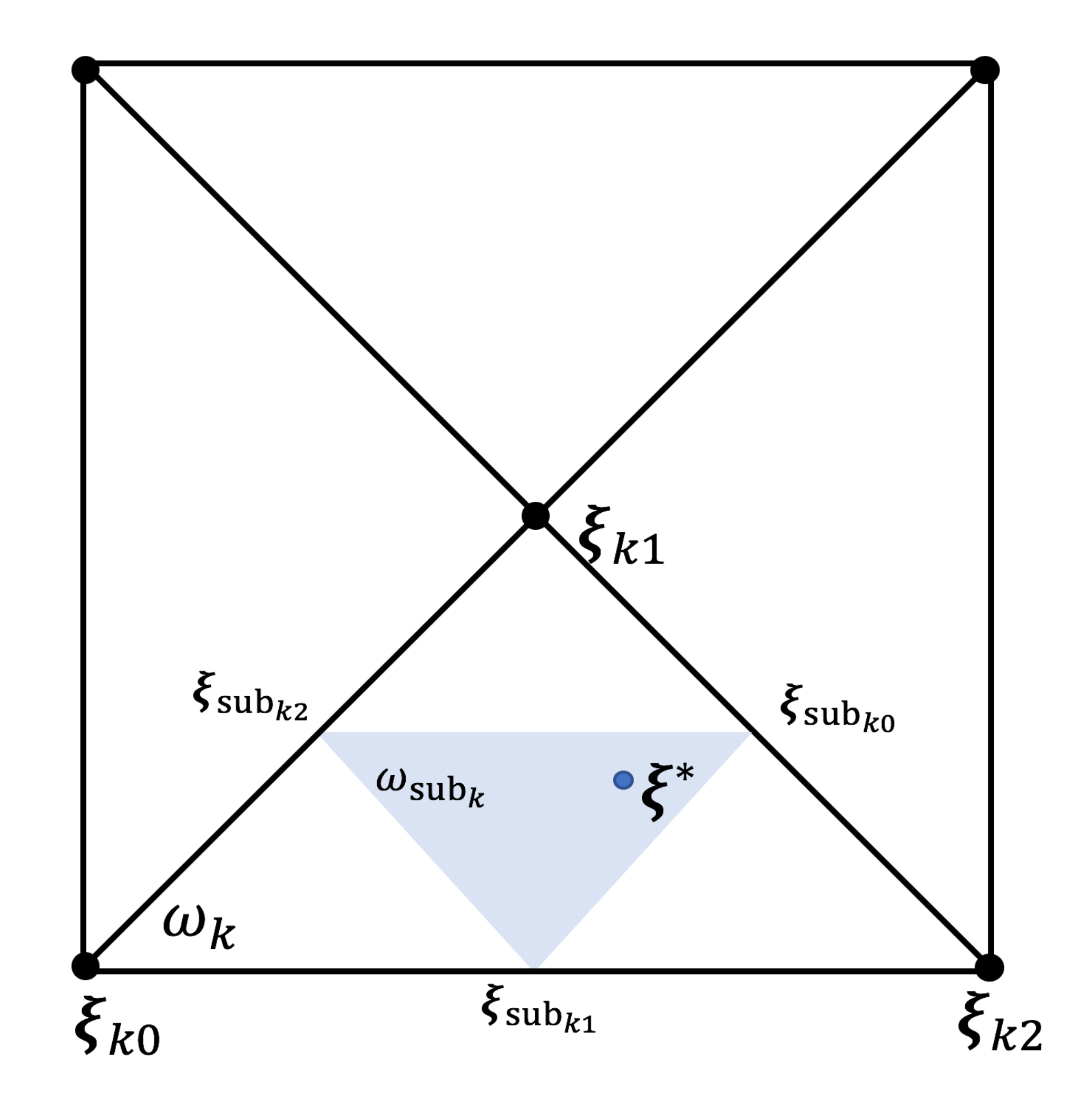}
	\caption{A sub-simplex (shaded area) in a two-dimensional case. The location of the random point (red dot) is bounded within the sub-simplex. }
	\label{subsimplex}
\end{figure}

Refinement of a selected simplex $\omega_{k}$ is performed by adding a new sample point $\boldsymbol{\xi}^\star$ in its interior.  The exact location of the newly added point is random and confined to a sub-simplex $\omega_{\text{sub}_{k}}$ of $\omega_{k}$,  in order to ensure sufficient spreading of samples (i.e. $\boldsymbol{\xi}^\star$ will not be close to the vertices of $\omega_{k}$). In the general case the $n_d+1$ vertices $\boldsymbol{\xi}_{\text{sub}_{j,l}}$ of $\omega_{\text{sub}_{k}}$ are defined as the centers of the faces of $\Omega_{k}$ 
\begin{equation}\label{eq:18}
\boldsymbol{\xi}_{\text{sub}_{k,l}} = \frac{1}{n_d}\sum_{\substack{l^{*}=0 \\  l^{*}\neq l}}^{n_d} \boldsymbol{\xi}_{k,l^*}, \quad l=0,\ldots,n_d
\end{equation}
Fig.\ \ref{subsimplex} depicts a sub-simplex for $n_d=2$. Upon refinement, the location of a new point (red dot in Fig.\ \ref{subsimplex}) is randomly selected in the interior of the sub-simplex \cite{SSC2010}.

According to Eq.\ \eqref{eq:18} new samples will fall only in the simplex interiors and the boundaries of the hypercube will never be sampled. This may lead to large aspect ratio elements that do not resolve features along the edges well. To avoid that, if a simplex selected for refinement is located at the boundary of the probability space then the new sample will be placed randomly along the boundary with probability 0.5. Otherwise, it is selected from within the sub-simplex.


%


\section{Proposed method}
\label{S:4}

\subsection*{Problem definition}
\noindent
Consider a physical system described by a set of governing differential equations and having stochastic input, system parameters, and/or boundary conditions described by $\boldsymbol{\xi}$. The full-field solution (e.g.\ displacement, strain, or stress fields) obtained from the solution of the governing equations (typically using a numerical approach, e.g.\ finite element method) is assembled in a matrix form $\textbf{F} \in \mathbb{R}^{n_f \times m_f}$ where  $n_f \times m_f=n$ is the number of degrees of freedom of the solution. In general, this matrix is assembled by arbitrarily reshaping the solution. For example, a typical FE code will produce a vector solution of dimension $n\times 1$ which can be reshaped into a matrix of arbitrary shape $n_f\times m_f=n$. Similarly, tensor-valued response (as may arise from time-varying matrix response) can be reassembled in matrix form by tensor unfolding. This reshaping should be performed with some care as it will affect both the dimension of the singular vectors and the rank $r$ of the matrix. In other words, it will affect the dimension of the Grassmann and Stiefel manifolds.

Let us define a reduced order form of $\textbf{F} \in \mathbb{R}^{n_f \times m_f}$  by performing SVD
\begin{equation}\label{eq:19}
\textbf{F} = \textbf{U}\boldsymbol{\Sigma}\boldsymbol{\Phi}^\intercal
\end{equation}
where $\boldsymbol{\Psi}\in \mathbb{R}^{n_f \times r}$ are the left-singular vectors, $\boldsymbol{\Sigma}\in \mathbb{R}^{r \times r}$ is a diagonal matrix of the singular values, and $\boldsymbol{\Phi}\in \mathbb{R}^{m_f \times r}$ are the right-singular vectors, all with rank $r$. Since matrices $\boldsymbol{\Psi}$ and $\boldsymbol{\Phi}$ are sets of orthonormal vectors of dimension $r$ in $\mathbb{R}^{n_f}$  and $\mathbb{R}^{m_f}$ respectively, we can consider them as points on the Stiefel manifolds $\mathcal{V}(r, n_f)$  and $\mathcal{V}(r, m_f)$, representing equivalent classes (points) on the Grassmann manifolds $\mathcal{G}(r,n_f)$ and $\mathcal{G}(r,m_f)$.


Notice that, because the system has stochastic properties, the solution $\textbf{F}$ is a random function of the parameters $\boldsymbol{\xi}$, (i.e. $\textbf{F} \approx \textbf{F}(\boldsymbol{\xi})$ it is a random field) which under Eq.\ \eqref{eq:19} defines $\boldsymbol{\Psi}(\boldsymbol{\xi}), \boldsymbol{\Sigma}(\boldsymbol{\xi})$ and $\boldsymbol{\Phi}(\boldsymbol{\xi})$. Thus, for every realization $i$ of the input parameter vector we have snapshots of $\boldsymbol{\Psi}_i = \boldsymbol{\Psi}(\boldsymbol{\xi}_i) \in \mathbb{R}^{n_f \times r_i}, \boldsymbol{\Sigma}_i =  \boldsymbol{\Sigma}(\boldsymbol{\xi}_i) \in \mathbb{R}^{r_i \times r_i}$ and $\boldsymbol{\Phi}_i =\boldsymbol{\Phi}(\boldsymbol{\xi}_i)\in \mathbb{R}^{m_f \times r_i}$0.  obtained from the SVD of $\textbf{F}_i = \textbf{F}(\boldsymbol{\xi}_i) \in \mathbb{R}^{n_f \times m_f}$.  

The  main idea of the proposed approach is to ``measure'' the difference in the behavior of the system using distances on the Grassmann manifold. Consider, for example, realizations $\boldsymbol{\xi}_i$ and $\boldsymbol{\xi}_j$ for the stochastic input to the system. The difference in the response of the system $\textbf{F}(\boldsymbol{\xi}_i)$ and $\textbf{F}(\boldsymbol{\xi}_j)$ is respresented by the distance between points $\mathcal{X}_i$ and $\mathcal{X}_j$ on the Grassmann manifold, represented by the orthonormal matrices $\boldsymbol{\Psi}_i$ and $\boldsymbol{\Psi}_j$ on the Stiefel manifold obtained from Eq.\ \eqref{eq:19}. These distances are used to identify regions of the space, $\Omega$, over which the behavior of the system changes rapidly. Samples are drawn more densely in these regions in order to facilitate an accurate approximation, $\tilde{\textbf{F}}$, of $\textbf{F}$ that is determined from interpolation on the Grassmann manifold and approximates the full solution without the need to solve the expensive governing equations. In this sense, the interpolation serves as a ``surrogate model'' on the full solution derived entirely from the samples.

\subsection*{Refinement strategy}

\begin{figure}[!ht]
	\centering\includegraphics[width=0.75\columnwidth]{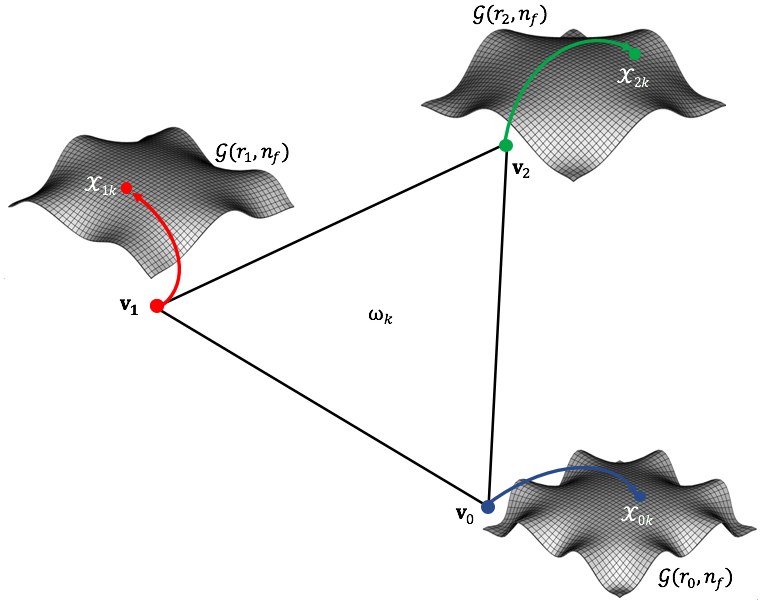}
	\caption{An element $\omega_k$ for a two-dimensional random input. Each vertex of $\omega_k$ corresponds to a point on a Grassmann manifold.}
	\label{Gr01}
\end{figure}

The method has a multi-element character as described in Section 3 since, as a product of the DT of the  probability space, every element $\omega_k$ is independent and can be fully defined by the location of its vertices (sample realizations) in the probability space and their corresponding points on $\mathcal{G}(r,n_f)$. For example,  if $n_d=2$ then simplex  $\omega_k$ is a triangle, defined by its three vertices  $\textbf{v}_{ik}\in \mathbb{R}^2$ and their  corresponding  Grassmann points $\mathcal{X}_{ik} \in \mathcal{G}(r_i, n_f)$, for $i=0,1,2$ (Fig.(\ref{Gr01})). 

At level $m$ of the adaptive procedure new points $\boldsymbol{\xi}^\star_i$ are added to the probability space and the current DT mesh is refined as described in Section 3. That is, from the set of $n_e$ elements $\Omega = \bigcup_{k=0}^{n_e-1}\omega_k$ at level $m$, a subset $\Omega_s \subseteq\Omega$ of $n_{ref}$ elements is selected for refinement according to an appropriate distance metric. The selection of a distance metric is described herein.

Two general strategies can be followed to compute the distance between points on the Grassman manifold that differ in how the dimension of the Grassmann manifold is identified (i.e.\ how the rank of $\boldsymbol{\Psi}_i$ is assigned). In the first approach a global rank $r_{\text{glob}}$ is provided a priori, which will be fixed for every SVD performed during the refinement procedure. In this case, there exists a unique Grasmmann manifold $\mathcal{G}(r_{\text{glob}},n_f)$ on which all the distances are calculated using e.g.\ the distances in Table 1.  The second strategy uses a tolerance for rank estimation of each snapshot $\mathbf{F}_i$ defined by $\text{tol} = \max({\Sigma}_i) \times n_f \times \varepsilon$, where $\max({\Sigma}_i)$ is the maximum singular value, $n_f$ is the dimension of the high-fidelity solution and $\varepsilon$ is the machine epsilon. Alternatively the user can prescribe a tolerance level for the SVD. By prescribing a tolerance level, $\text{tol}$, the solution at each point possesses a unique rank and therefore lies on its own Grassmann manifold requiring the computation of distances on the doubly infinite Grassmann manifold (Table 2). A detailed investigation of these strategies will be considered in the numerical example.

The refinement metric adopted in this work is the total Grassmann distance along the perimeter of each element, expressed as the sum of all pairwise Grassmann distances between its vertices $\boldsymbol{\Psi}_{ik}$:
\begin{equation}\label{eq:20}
D_k = \sum_{\substack{i, j=0 \\ i\neq j}}^{n_d} \delta(\boldsymbol{\Psi}_{ik},\boldsymbol{\Psi}_{jk})
\end{equation}
where $n_d$ is the number of random variables used to describe the uncertainties. Again, the definition of $\delta$ will depend upon the relative dimensions of the Grassmann manifolds of the vertices. For strategy one, where the dimension $r_{\text{glob}}$ is fixed, we have $\delta=d_{\mathcal{G}(r_{\text{glob}},n_f)}(\boldsymbol{\Psi}_{ik},\boldsymbol{\Psi}_{jk})$. Meanwhile, for strategy two where a tolerance is specified and the dimension varies among the vertices, $\delta=d_{\mathcal{G}(\infty,\infty)}(\boldsymbol{\Psi}_{ik},\boldsymbol{\Psi}_{jk})$. At each level the elements with the highest values of $D_k$ are selected (i.e.\ $\omega_k:D_k>D_{th}$) for refinement. This has the effect of concentrating samples in regions of the probability space where the distances between different snapshots of the solution are large, i.e the model exhibits a significant change in its behavior.

\subsection*{Interpolation on the manifold}
\label{Interpolation}

A key ingredient of the proposed method is the ability to interpolate on the Grassmann manifold to ``predict'' the full system solution at previously unsampled points without solving the governing equations. This interpolation is performed on each DT element using the method proposed by Amsallem and Farhat \cite{Farhat2008}. The major difference between the interpolation performed herein and that performed in \cite{Farhat2008} is that the authors of \cite{Farhat2008} use interpolation to identify a reduced basis upon which to project the governing equations (thus making the solution less expensive) whereas we perform a full interpolation of the solution - thereby eliminating the reduced order model. Details of this method follow.


Recall the tangent space at point $\mathcal{X} \in \mathcal{G}(r,n_f)$ (Eq.\ \eqref{eq:3}) is a vector space, and thus, it is a flat space on which interpolation between two bases can be performed in a standard way. Consider the ${n_d}+1$ vertices of element $\omega_k$ generated from realizations of the $n_d$-dimensional random vector $\boldsymbol{\xi}_i$  with $i = 0, \ldots, n_d+1$ and with corresponding points on the Grassmann manifold $\mathcal{X}_i\in \mathcal{G}(r_i,n_f)$, spanned by the columns of $\boldsymbol{\Psi}_i \in \mathbb{R}^{r_i\times n_f}$, respectively. In general, $r_i\ne r_j$ for $i \neq j$ and therefore $\mathcal{X}_i$ and $\mathcal{X}_j$ lie on different manifolds.  To perform interpolation these points must lie on the same manifold. To achieve this we retain only those basis vectors that form the sets $\mathcal{X}_0,\ldots,\mathcal{X}_{n_s} \in \mathcal{G}(\min(r_i;~i=1,\dots,n_d+1),n_f)$. Selecting arbitrarily $\boldsymbol{\Psi}_0$ as the origin Eq.\ \eqref{eq:7} becomes
\begin{equation}\label{eq:21}
\mathbf{M}_i=(\textbf{I}-\boldsymbol{\Psi}_0\boldsymbol{\Psi}_0^\intercal)\boldsymbol{\Psi}_i(\boldsymbol{\Psi}_0^\intercal \boldsymbol{\Psi}_i)^{-1} = \mathbf{U}_i \textbf{S}_i\mathbf{V}_i^\intercal
\end{equation}
where $\mathbf{U}_i,  \textbf{S}_i$ and $\mathbf{V}_i$ are the left eigenvectors, singular values and right-singular vectors, respectively, of matrix $\textbf{M}_i$.  Thus, the matrix $\boldsymbol{\Gamma}_i$ representing point $\mathcal{H}_i$ on $\mathcal{T}_{\mathcal{X}_0}$ is given by
\begin{equation}\label{eq:22}
\boldsymbol{\Gamma}_i = \mathbf{U}_i\tan^{-1}(\textbf{S}_i)\mathbf{V}_i^\intercal
\end{equation}

For a new realization of the random vector $\boldsymbol{\xi}^\star$, 
the interpolated matrix $\tilde{\boldsymbol{\Gamma}} \in \mathcal{T}_{\mathcal{X}_0}$ is estimated by

\begin{equation}\label{eq:23}
\tilde{\boldsymbol{\Gamma}} = \frac{1}{V_{\omega_k}}\sum_{i = 0}^{n_d+1} V_{\omega_{\text{sub}_k}^i}\boldsymbol{\Gamma}_i 
\end{equation}
where $V_{\omega_k}$ is the volume of element $\omega_k$ estimated as
\begin{equation}
V_{\omega_k}=\dfrac{1}{n_d!}\det\bigg({\left[\boldsymbol{\xi}_{l_{k,1}}-\boldsymbol{\xi}_{l_{k,0}}, \boldsymbol{\xi}_{l_{k,2}}-\boldsymbol{\xi}_{l_{k,0}}, \dots, \boldsymbol{\xi}_{l_{k,n_d+1}}-\boldsymbol{\xi}_{l_{k,0}}\right]\bigg)}
\end{equation}
and the subscript $l_{k,j}$ denotes the global point $\boldsymbol{\xi}_l$ serving as the $j^{th}$ vertex of element $\omega_k$. Similarly, $V_{\omega_{\text{sub}_k}^i}$ is the volume of the sub-simplex $\omega_{\text{sub}_k}^i$ in $\omega_k$ defined by replacing vertex $\boldsymbol{\xi}_i$ with the interpolation point $\boldsymbol{\xi}^\star$.  Other, more advanced interpolation methods can be considered using, for example, polynomial interpolating functions but we have found this linear barycentric interpolation sufficient for our purposes.

Combining Eq.\ \eqref{eq:22} and Eq.\ \eqref{eq:23} and keeping in mind that $\boldsymbol{\Gamma}_0 = \textbf{0}$ since $\mathcal{X}_0$ is the origin of the tangent space we get
\begin{equation}\label{eq:24}
\tilde{\boldsymbol{\Gamma}} = \frac{1}{V_{\omega_k}}\sum_{i = 1}^{n_d+1}\Big(\mathbf{U}_i\left[V_{\omega_{\text{sub}_k}^i}\tan^{-1}(\textbf{S}_i)\right]\mathbf{V}_i^\intercal\Big)
\end{equation}





After estimating $\tilde{\boldsymbol{\Gamma}}$, its projection on $\mathcal{G}(\min(r_0,\ldots,r_{n_s-1}),n_f)$ is given by the exponential mapping as \cite{Farhat2008}
\begin{equation}\label{eq:25}
\tilde{\boldsymbol{\Psi}} =  \left(\boldsymbol{\Psi}_0\tilde{\mathbf{V}}\cos(\tilde{\mathbf{S}})+\tilde{\mathbf{U}}\sin(\tilde{\mathbf{S}}) \right)\tilde{\mathbf{V}}^\intercal
\end{equation}
where $\tilde{\boldsymbol{\Gamma}}=\tilde{\boldsymbol{\Phi}}\tilde{\mathbf{S}}\tilde{\boldsymbol{\Psi}}^\intercal$ and whose columns span  $\mathcal{X}^\star \in \mathcal{G}(\min(r_0,\ldots, r_{n_s-1}), n_f)$.



Eq.\ \eqref{eq:25} identifies a basis for the solution at $\boldsymbol{\xi}^\star$ but does not uniquely define the solution. To estimate the full solution at this point also requires both the right singular vectors and the singular values. Fortunately, the right singular vectors similarly define points on $ \mathcal{G}(\min(r_0,\ldots, r_{n_s-1}), n_f)$, which may be similarly interpolated to yield $\tilde{\boldsymbol{\Phi}}$. Meanwhile, under conditions of sufficiently small distances $\delta(\boldsymbol{\Psi}_i,\boldsymbol{\Psi}_j)$ and $\delta(\boldsymbol{\Phi}_i,\boldsymbol{\Phi}_j)$, we assume that the singular values can be directly interpolated to estimate $\tilde{\boldsymbol{\Sigma}}$ thus completing the full solution interpolation as $\tilde{\textbf{F}} = \tilde{\boldsymbol{\Psi}}\tilde{\boldsymbol{\Sigma}}\tilde{\textbf{V}}$. In the subsequent example, this is shown to produce a very accurate estimate when interpolating between points that are sufficiently close on the Grassmann.

\subsection*{Stopping criterion}

To control the refinement procedure, an error estimate is required which establishes a metric of accuracy for the interpolation and provides a stopping criterion. The error estimates adopted in this work is based on an average principal angle, $\tilde{\theta}$, between the interpolated solution and the true solution estimated at each newly added point $\boldsymbol{\xi}^\star$ inside the sub-simplex $\omega_{\text{sub}_k}$ of element $\omega_k$.  This average principal angle is defined as 
\begin{equation}\label{eq:26}
\tilde{\theta}(\boldsymbol{\xi}^\star) = \sqrt{\frac{\delta(\boldsymbol{\Psi}_k(\boldsymbol{\xi}^\star),\tilde{\boldsymbol{\Psi}}_k(\boldsymbol{\xi}^\star))}{\min(r,\tilde{r}_k)}}
\end{equation}
where $\tilde{\boldsymbol{\Psi}}_k(\boldsymbol{\xi}^\star)$ and $\boldsymbol{\Psi}_k(\boldsymbol{\xi}^\star)$ are the approximate basis (with rank $\tilde{r}_k$) estimated by interpolation on the manifold and the actual basis (with rank $r$) corresponding to the model evaluation, respectively. Again, $\delta$ may be either the classical or doubly infinite Grassmann distance depending on the Grassmann dimensions. 

In refinement level $l+1$, if an element $\omega_k$ has at least $n_d$ out of $n_d+1$ vertices $\mathcal{X}_{ik}$ for which $\tilde{\theta}$ was less than a threshold value in refinement level $l$, i.e.
\begin{equation}\label{eq:27}
\tilde{\theta} \leq \theta_{\text{ref}}
\end{equation}
then this element is considered sufficiently converged and is no longer refined. That is, for all subsequent refinement levels, $l+1,l+2,\dots$ elements with $n_d$ out of $n_d+1$ vertices are excluded from further refinement. The threshold principal angle $\theta_{\text{ref}}$ is bounded on $0\leq \theta_{\text{ref}}\leq \pi/2$ such that the larger the value of $\theta_{\text{ref}}$ the lower quality of the approximation, in general. 

The overall procedure is terminated when all elements are considered sufficiently converged. In this work, we have used the value $\theta_{\text{ref}}=\pi/15$ which gives accurate results.

A brief description of the proposed method can be found in Algorithm 1.
\begin{algorithm}[!ht]
	\caption{Proposed adaptive method using Grassmann variations }\label{euclid}
	\begin{algorithmic}[1]
		\Procedure{given $\boldsymbol{\xi}_1, \boldsymbol{\xi}_2,\ldots,\boldsymbol{\xi}_{n_s}$}{} 
        \BState $\{\textbf{F}_1, \textbf{F}_2,\ldots,\textbf{F}_{n_s}\} \leftarrow \{\boldsymbol{\xi}_1, \boldsymbol{\xi}_2,\ldots,\boldsymbol{\xi}_{n_s}\}$
        \For{$i =1, \ldots, n_s$}
        \State $\textbf{F}_i = \boldsymbol{\Psi}_i\boldsymbol{\Sigma}_i\boldsymbol{\Phi}_i^\intercal$ (thin SVD)
        \EndFor
         \BState $n_e \leftarrow n_{e}=O(n_{s}^{n_{d}/ 2})$
         \For{$k =0, \ldots, n_e-1$}
        \State $\delta(\boldsymbol{\Psi}_{ik},\boldsymbol{\Psi}_{jk}) \leftarrow$ Eq.\eqref{eq:17}  where $i,j \in [0,\ldots,n_d]$
        \State $D_k \leftarrow$ Eq.\eqref{eq:20} 
        \EndFor
       \BState $n_c \leftarrow$ elements with at most $n_d$ out of $n_d+1$ vertices having $\tilde{\theta} \leq \theta_{\text{ref}}$
       \BState From $n_c$ candidates, find $n_{ref}$ elements for which $D_k \geq D_{th}$
 	   \BState  $i \leftarrow 0$
       \While{ $i \leq n_{ref}-1$} 
       \State Generate random sample $\boldsymbol{\xi}^\star$ in $\omega_{\text{sub}_i}$ (or on boundary)
       \State $\tilde{\theta}(\boldsymbol{\xi}^\star) \leftarrow$ Eq.\eqref{eq:26} 
       \State $i \leftarrow i+1$
       \BState \textbf{Repeat} from Step 1 until $n_c=0$
       \EndWhile

		\EndProcedure
	\end{algorithmic}
\end{algorithm}

\subsection*{Extension to problems with high-dimensional input}
The proposed method is limited by the Delaunay triangulation to problems of modest stochastic dimension ($n_d$ less than about 6). Extension to higher-dimensional problems will require adaptation of the proposed methodology to employ a high-dimensional model representation (HDMR) or similar decomposition of the model as performed in \cite{Ma2010,Edeling16} or to abandon the multi-element discretization of the input. Investigations of these approaches are ongoing and will be considered in a future work.

\section{Application to material modeling}
\label{S:5}

\noindent
The model used to demonstrate the applicability and the efficiency of the proposed methodology is an elasto-viscoplastic material model for a bulk metallic glass (BMG) based on the shear transformation zone (STZ) theory of amorphous plasticity \cite{Falk1998, Langer2008}, using  a quasi-static integration method  for the numerical solution \cite{Rycroft2015}.

\subsection{Elasto-plastic material model}
\label{S:5.1}

\noindent
For an elastoplastic material the rate-of-deformation tensor $\textbf{D} = (\nabla \textbf{v}+ \nabla \textbf{v}^\intercal)/2$, where $\textbf{v}$ is the velocity, can be expressed as the sum 
$\textbf{D}= \textbf{D}^{el}+\textbf{D}^{pl}$ and subsequently the linear elastic constitutive relation is
\begin{equation}\label{eq:28}
\frac{\mathcal{D}\boldsymbol{\sigma}}{\mathcal{D}t} = \textbf{C}:\textbf{D}^{el} = \textbf{C}:(\textbf{D}-\mathbf{D}^{pl})
\end{equation}
where  $\textbf{C}$ is the fourth-order stiffness tensor, (assumed to be isotropic, homogeneous and time-independent), $\boldsymbol\sigma$ is the Cauchy stress tensor and $\mathbf{D}^{el}$ and $\mathbf{D}^{pl}$ are the elastic and plastic rate-of-deformation tensors, respectively. Under the assumption that the elastic deformation is relatively small and taking into account the translation and the rotation of the material, the time-evolution of the stress is given by the Jaumann objective stress rate $\mathcal{D}\boldsymbol{\sigma}/\mathcal{D}t = \mbox{d}\boldsymbol{\sigma}/\mbox{d}t + \boldsymbol{\sigma}\cdot \boldsymbol{\omega} -\boldsymbol{\omega} \cdot \boldsymbol{\sigma} $, where $\mbox{d}/\mbox{d}t$ is the advective derivative that takes into account material movement, i.e. $\mbox{d}/\mbox{d}t = \partial/\partial t + (\mathbf{v}\cdot \nabla)$ and $\boldsymbol{\omega}=(\nabla \textbf{v}- \nabla \textbf{v}^\intercal)/2$ is the spin. By considering force balance, the velocity satisfies

\begin{equation}\label{eq:29}
\rho\frac{\mbox{d}\textbf{v}}{\mbox{d}t} = \nabla \cdot \boldsymbol\sigma
\end{equation}
where $\rho$ is the density of the material. Using a finite-difference numerical solution with explicit time integration, one can solve the system of the hyperbolic differential equations (\ref{eq:28}) and (\ref{eq:29}). However, such a solution scheme can be computationally expensive since prohibitively small timesteps and are required in order to satisfy the Courant-Friedrichs-Lewy (CFL) condition \cite{Heath}.

Because the plastic deformation occurs on a much longer time-scale than elastic wave propagation Rycroft et al.\ \cite{Rycroft2015} argue that the acceleration term $\mbox{d}\textbf{v}/\mbox{d}t$ becomes negligible yielding a divergence free stress field

\begin{equation}\label{eq:30}
\nabla \cdot \boldsymbol\sigma = \mathbf{0}
\end{equation}
which implies that the forces remain in quasi-static equilibrium. The authors of \cite{Rycroft2015} recognized a strong similarity between this elasto-plastic model under the quasi-static condition and the incompressible Navier-Stokes equations and they developed an Eulerian numerical approach to solve this quasi-elasto-plastic system of equations which is adopted herein.

\subsection{Shear transformation zone theory of plasticity}
\label{S:5.2}

\noindent
In the shear transformation zone theory \cite{Falk1998}, plastic deformation occurs through highly localized rearrangements in the configuration of clusters of atoms called shear transformation zones (STZs), which change from one stable configuration to another according to a stress-biased thermal activation process. Consider a bulk metallic glass (BMG) - an amorphous metal that undergoes a glass transition - at a given temperature $T$ below the glass transition. The STZ theory assumes that the material is populated with STZs whose density is characterized in terms of an effective disorder temperature defined as:
\begin{equation}\label{eq:31}
\chi=\dfrac{\partial U_c}{\partial S_c}
\end{equation}
where $U_c$ and $S_c$ are the potential energy and the entropy of the configurational degrees of freedom\footnote{The effective disorder temperature is set on the concept that time-scale separation for the vibrational and configurational degrees of freedom is possible.}, respectively. In terms of the effective temperature $\chi$, a flow rule for the plastic rate-of-deformation tensor has been formulated:
\begin{equation}\label{eq:32}
\mathbf{D}^{pl}=\dfrac{1}{\tau_0}\exp\left\{-\left(\dfrac{-e_z}{k_B\chi}+\dfrac{\Delta}{k_B T}\right)\right\}\cosh\left(\dfrac{\tilde{\Omega}\epsilon_0\bar{\sigma}}{k_BT}\right)\left(1-\dfrac{\sigma_y}{\bar{\sigma}}\right)
\end{equation}
where $\sigma_y$ is the material's yield stress, $T$ is the classical thermodynamic temperature, $\bar{\sigma}=|\boldsymbol{\sigma}_0|$ is the magnitude of the deviatoric stress $\boldsymbol{\sigma}_0=\boldsymbol{\sigma}-\nicefrac{1}{3}\boldsymbol{1}\text{tr}\boldsymbol{\sigma}$, and $k_B=1.38\times 10^{-23}$ Joules/Kelvin is the Boltzmann constant. The rest of the model parameters $\tau_0$, $\epsilon_0$, $\Delta$, $\tilde{\Omega}$, and $e_z$ are explained in Table \ref{tab:plasticity} where their values are set in accordance with values reported in the literature \cite{Langer2008,Yu2010,Rycroft2012,Rycroft2015}. Similarly the elasticity parameters of the model are set in Table \ref{tab:elasticity}. The effective temperature evolves according to a heat equation as
\begin{equation}\label{eq:33}
c_0 \dot{\chi}=\dfrac{1}{\sigma_y}(\mathbf{D}^{pl}:\boldsymbol{\sigma}_0)(\chi_\infty-\chi)+\nabla\cdot D_\chi\nabla\chi
\end{equation}
where $c_0$ plays the role of a specific heat, $\chi_\infty$ is a saturation temperature that represents the state of maximum disorder and $D_\chi=l_d^2\sqrt{\mathbf{D}^{pl}:\mathbf{D}^{pl}}$ is a rate dependent diffusivity; parameter $l_d$ sets a length scale for plastic deformation mediated diffusion. 
\begin{table}[!ht]
	\centering
	\caption{Parameters for the STZ plasticity model for Vitreloy 1 bulk metallic glass. Note that $k_B=1.38\times 10^{-23}\hspace{3pt} \mbox{JK}^{-1}$ is the Boltzmann constant.}
	\begin{tabular}{ccl}
	\hline
	Parameter & Value & Description \\
	\hline \hline
	$\sigma_y$ & 0.85 GPa & Yield stress \\
	$\tau_0$ & 1e-13 s & Molecular vibration timescale \\
	$\epsilon_0$ & 0.3 & Typical local strain at STZ transition \\
	$\Delta$ & $8000k_B ~\textrm{K}$  & Typical activation barrier \\
	$\tilde{\Omega}$ & 300 \AA & Typical activation volume \\
	$T$ & 400 K & Classical thermodynamic temperature \\
	$\chi_\infty$ & 900 K & Steady-state effective temperature \\
	$e_z$ & $21,000k_B$ K & STZ formation energy \\
	$c_0$ & 0.4 & Specific heat-like fraction of plastic  \\
    &            & work that increases $\chi$ \\
	$l_d$ & 0.06 cm & Length-scale for rate-dependent diffusivity\\
	\hline
	\end{tabular}
	\label{tab:plasticity}
\end{table}

\begin{table}[!ht]
	\centering
	\caption{Elasticity parameters for Vitreloy 1 bulk metallic glass.}
	\begin{tabular}{ccl}
	\hline
	Parameter & Value & Description \\
	\hline \hline
	$E$ & 101 GPa & Young's modulus \\
	$\nu$ & 0.35 & Poisson's ratio \\
	$K$ & 122 GPa & Bulk modulus \\
	$\mu$ & 37.4 GPa & Shear modulus \\
	$\rho$ & 6125 $\text{kg m}^{-3}$ & Density \\
	$c_s=\sqrt{\mu/\rho}$ & 2.47 $\text{km s}^{-1}$ & Shear wave speed \\
	\hline
	\end{tabular}
	\label{tab:elasticity}
\end{table}

In this work, the STZ plasticity model is simulated using the quasi-static numerical method described above and detailed in \cite{Rycroft2015}, for a two-dimensional square specimen of Vitreloy 1 BMG with dimensions 2 cm $\times$ 2 cm using a grid spacing $\Delta x=0.01$ cm subjected to simple shear. Hence, the computational domain consists of a square grid of $201 \times 201$ integration points (i.e.\ the solution matrix $\mathbf{F}$ has dimension $201 \times 201$ and reshaping is not necessary). Simple shear is imposed by moving the top and bottom boundaries of the system with a fixed velocity $U_B=247.13$ cm/s for a duration of $t^*=8.09e-4$ s yielding a total displacement of 2 mm.

\subsection{Uncertainty quantification for full strain field response of a bulk metallic glass}

\noindent
Past studies indicate that the STZ plasticity model is highly sensitive to the initial effective temperature field, which can lead to either highly localized shear deformation or relatively homogeneous plastic deformation \cite{Manning2007,Manning2009,Shields}. Since little data exists to suggest appropriate parameters for the initial $\chi$ field \cite{Hinkle}, a high degree of uncertainty exists in the expected material behavior. More specifically, the initial $\chi$ field at $t=0$ is coarse-grained from atomistic simulations \cite{Hinkle}, 
but the parameters of the coarse-graining are uncertain with two stochastic dimensions: the mean effective temperature $\mu_\chi$ and the coefficient of variation of the effective temperature $c_\chi$, i.e. $\textbf{s} = [\mu_\chi, c_\chi]$, with probability distributions given in Table \ref{tab:estochasticity}.  This coarse-graining means that the $\chi$ value at each grid-point comes from lower-scale information 
by convolving a fixed grid of independent normalized atomic potential energies with a Gaussian kernel having fixed standard deviation $l_\chi$ and scaling the field according to $\mu_\chi$ and $\sigma_\chi = c_\chi\cdot \mu_\chi $. Note that the length-scale $l_{\chi}$ may also be uncertain but is treated here as fixed for illustration purposes.
\begin{table}[!ht]
	\centering
	\caption{Random field parameters and their probability distributions.}
	\begin{tabular}{ccll}
	\hline
	Parameter & Distribution & Value & Description\\
	\hline \hline
	$\mu_\chi$ & Uniform & [500,700] & Mean effective temperature\\
	$c_\chi$ & Uniform & [0,0.1] & Coefficient of variation \\
    $l_\chi$  & - & 10 & Correlation length-scale \\
	\hline
	\end{tabular}
	\label{tab:estochasticity}
\end{table}
 
Since our stochastic problem is defined in a two-dimensional space $(n_d = 2)$, the initial DT grid  consists of four sample points located at the corners and one in the middle of the square probability space $C^2=[0,1]\times[0,1]$ (i.e. $\boldsymbol{\xi}_0= (0,0), \boldsymbol{\xi}_1 = (0,1), \boldsymbol{\xi}_2 = (1,0), \boldsymbol{\xi}_3 = (1,1), \boldsymbol{\xi}_4 = (0.5,0.5)$). To map the variables between the parameter space $S$ and the probability space $\Omega$ we use the following simple linear mapping: $\mu_\chi^{(i)}=500+200{\xi}_{i1}$, $c_\chi^{(i)}=0.1{\xi}_{i2}$.



\noindent

\begin{figure}[!ht]
	\centering\includegraphics[width=1.\columnwidth]{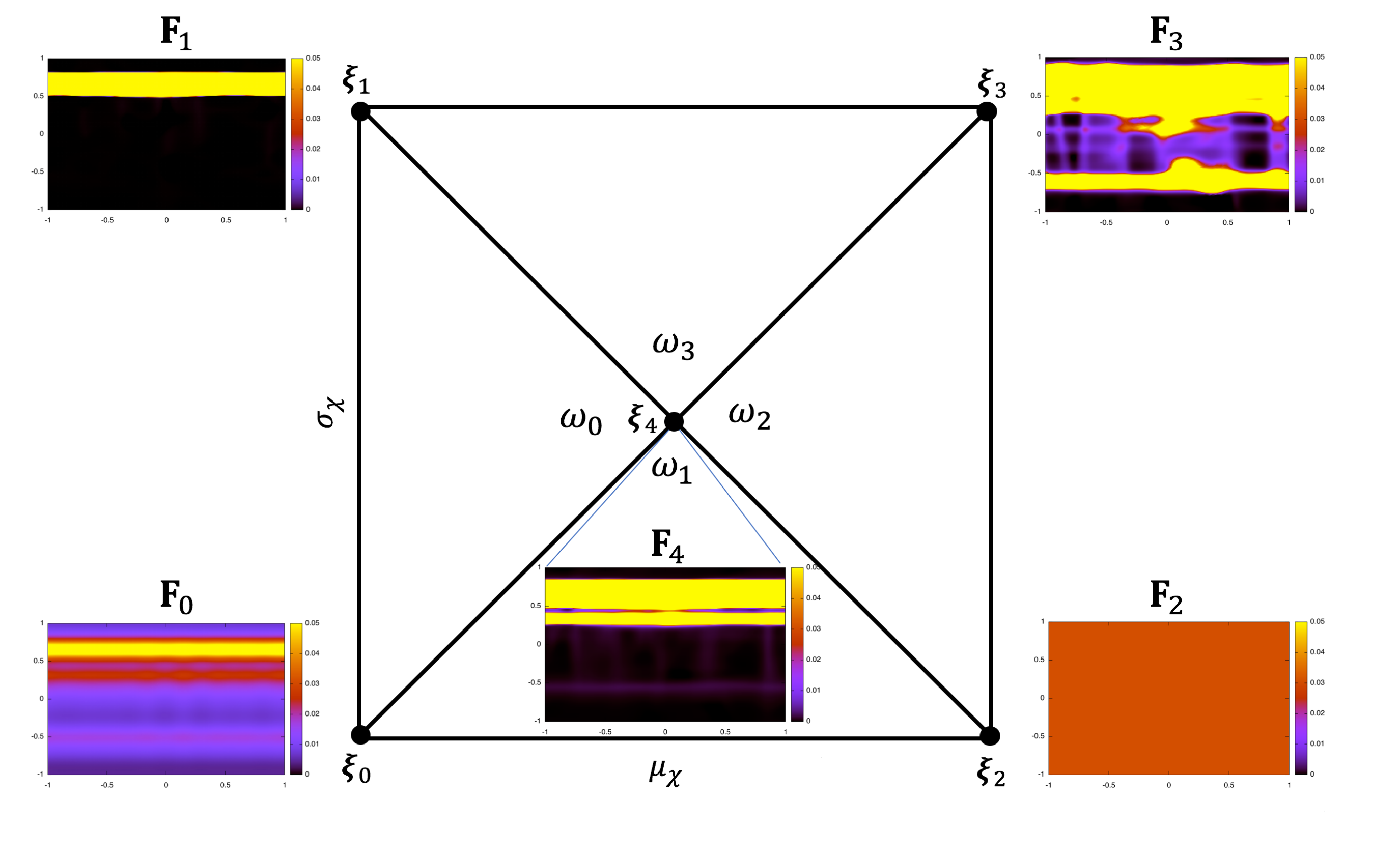}
	\caption{Delaunay triangulation of the initial $2^{n_d=2}+ 1= 5$ samples. Four triangular elements partition the probability space, each one having as vertices snapshots of the plastic strain field $\textbf{F}_i$.}
	\label{fig:DT0}
\end{figure}

At each of the initial sample points $\boldsymbol{\xi}_{i}$, the STZ equations are solved and the full plastic strain field retained. In keeping with previous notation, this ($201\times201$ dimensional) response field is denoted by $\textbf{F}_i$. To visualize the location in the probability space of these snapshots, Fig.\ \ref{fig:DT0} depicts the initial DT represented by four simplices $\omega_k$, $k=0,\ldots,3$ which discretize $\Omega=C^2$ into four elements, together with snapshots of the solution $\textbf{F}_i$ corresponding to each realization of the random vector $\mathbf{s}_i = [\mu_\chi^{(i)}, c_\chi^{(i)}]$, $i=0,\ldots,4$ resulting from transformation of $\boldsymbol{\xi}_i$. 
In the proposed methodology we use the metrics on $\mathcal{G}(\infty,\infty)$ given in Table \ref{tab:2} to compare the solutions on the vertices of each element. For example, simplex $\omega_0$ has vertices $\boldsymbol{\xi}_0$,  $\boldsymbol{\xi}_1$ and $\boldsymbol{\xi}_4$. To calculate the perimeter of $\omega_0$ on the $\mathcal{G}(\infty, \infty)$ we need to calculate the pairwise distances $\delta(\boldsymbol{\Psi}_0,\boldsymbol{\Psi}_1)$, $\delta(\boldsymbol{\Psi}_0,\boldsymbol{\Psi}_4)$ and $\delta(\boldsymbol{\Psi}_1,\boldsymbol{\Psi}_4)$, which represent the pairwise differences between the plastic strain fields, i.e. $\textbf{F}_{0}-\textbf{F}_{1}$, $\textbf{F}_{0}-\textbf{F}_{4}$ and $\textbf{F}_{1}-\textbf{F}_{4}$. Recall that $\boldsymbol{\Psi}_i$ are estimated from $\textbf{F}_{i}$ using the procedure described in Section \ref{S:4}. 

Table \ref{tab:3} shows the distances computed for each of the elements in Figure \ref{fig:DT0} according to three different metrics on $\mathcal{G}(\infty,\infty)$. Notice a slight change in the indexing of $\boldsymbol{\Psi}$ from a global indexing $\boldsymbol{\Psi}_i$ where $i$ denotes the vertex to a local indexing on element $k$, $\boldsymbol{\Psi}_{jk}$ where $j$ indexes on the vertex number local to the element (i.e. $j=0,1,2$) such that e.g. $\boldsymbol{\Psi}_{01}$ denotes vertex 0 on element $\omega_1$. The overall distance metric $D_k$ from Eq.\ \eqref{eq:20} is calculated for every element $k=0,\ldots, n_e-1$, and the elements with value larger than a threshold, i.e. $D_k \geq D_{th}$  are selected for refinement. We set the threshold $D_{th}$ as the $\alpha-$th percentile of the vector $\{D_0, \ldots, D_{n_e-1}\}$. Here, the value of $D_{th}$ corresponds to $\alpha=0.80$ such that, at each iteration level, 20\% of the elements are refined.
\begin{table}[!ht]
	\captionof{table}{Distances computed using different metrics on $\mathcal{G}(\infty,\infty)$ for each of the simplex elements from the DT of Fig.\ \ref{fig:DT0}} 
    \resizebox{\columnwidth}{!}{
	\begin{tabular}{ c c c c c c}   \toprule
	\textbf{Metric}	&\textbf{Simplex} & $\delta(\boldsymbol{\Psi}_{0k},\boldsymbol{\Psi}_{1k})$& $\delta(\boldsymbol{\Psi}_{0k},\boldsymbol{\Psi}_{2k})$ & $\delta(\boldsymbol{\Psi}_{1k},\boldsymbol{\Psi}_{2k})$ & $D_k$ (Eq.\ \eqref{eq:20}) \\\midrule
\multirow{ 4}{*}{\textbf{Grassmann}}&$\omega_0$	&	8.46 & 7.32 & 10.79 & 26.59 \\ 
	&	$\omega_1$         & 10.79 & 4.06 & 11.04 & 25.90      \\ 
	&	$\omega_2$      & 11.04 & 18.65& 15.86  & $\mathbf{45.55}$   \\ 
     &   $\omega_3$      & 15.86 & 7.32 & 17.40  & 40.59   \\ 
		\hline
        \hline
        
 \multirow{ 4}{*}{\textbf{Chordal}}	&	$\omega_0$         & 5.46 & 4.72 & 6.94 &  17.13                                  \\ 
	&	$\omega_1$         & 6.94 & 2.70 & 7.09 & 16.73      \\ 
	&	$\omega_2$      & 7.09 & 11.91 & 10.29  & $\mathbf{29.30}$    \\ 
     &   $\omega_3$      & 10.29 & 4.72& 11.21  & 26.23   \\ 
		\hline       
        \hline
   \multirow{ 4}{*}{\textbf{Procrustes}}	&	$\omega_0$         & 5.40 & 4.67 & 6.88 &    16.96                                \\ 
	&	$\omega_1$         & 6.88 & 2.60 & 7.04 & 16.53      \\ 
	&	$\omega_2$      & 7.04 & 11.87& 10.13  & $\mathbf{29.05}$   \\ 
     &   $\omega_3$      & 10.13 & 4.67& 11.10 & 25.91    
		\\\bottomrule    
		\label{tab:3}
	\end{tabular}}
\end{table}



\begin{figure}[!ht]
	\centering\includegraphics[width=1.0\columnwidth]{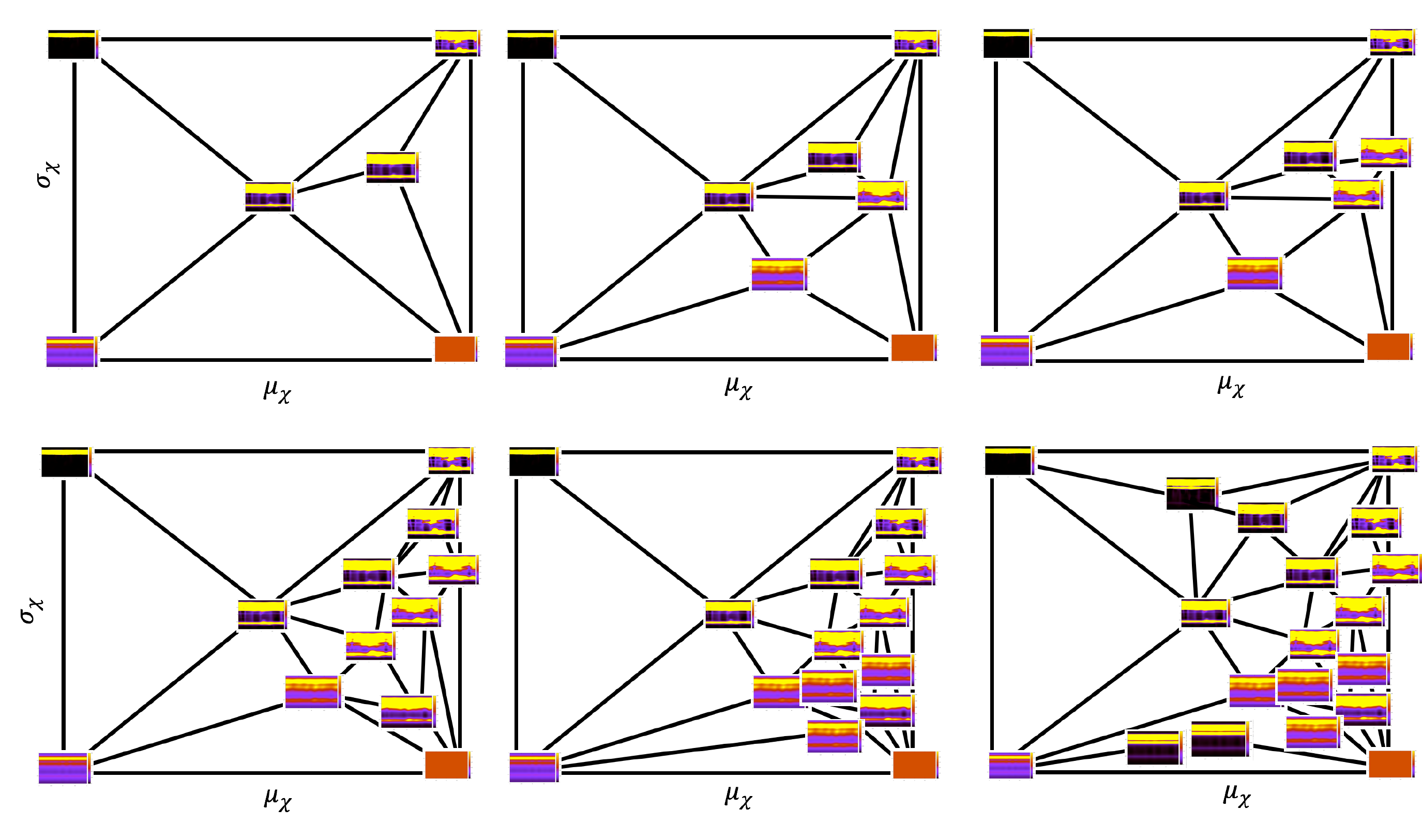}
	\caption{Evolution of the adaptive sampling showing the DT along with the plastic strain field at each of the vertices; The algorithm is capable of concentrating samples in areas of transition between homogeneous and localized plasticity from early stages. }
	\label{fig:DT}
\end{figure}

From Table \ref{tab:3} and Figure \ref{fig:DT0}, it is clear that element $\omega_2$ has the largest total distance - regardless of metric. Hence it is refined first.  Figure \ref{fig:DT} shows the DT refinement through six levels along with the corresponding plastic strain field at each vertex.  From  the initial levels of refinement, the algorithm identifies regions where the material behavior transitions from homogeneous plastic deformation (e.g.\ $\boldsymbol{\xi}_2=(1,0)$) to highly localized plastic strains in the form of a single dominant shear band (e.g.\ $\boldsymbol{\xi}_1=(0,1)$) or multiple shear bands (e.g.\ $\boldsymbol{\xi}_3=(1,1)$).
These transition regions are populated with additional samples, which aids in understanding the material response and improves the quality of approximations generated by manifold interpolations.

\begin{figure}[!ht]
	\centering\includegraphics[width=0.8\columnwidth]{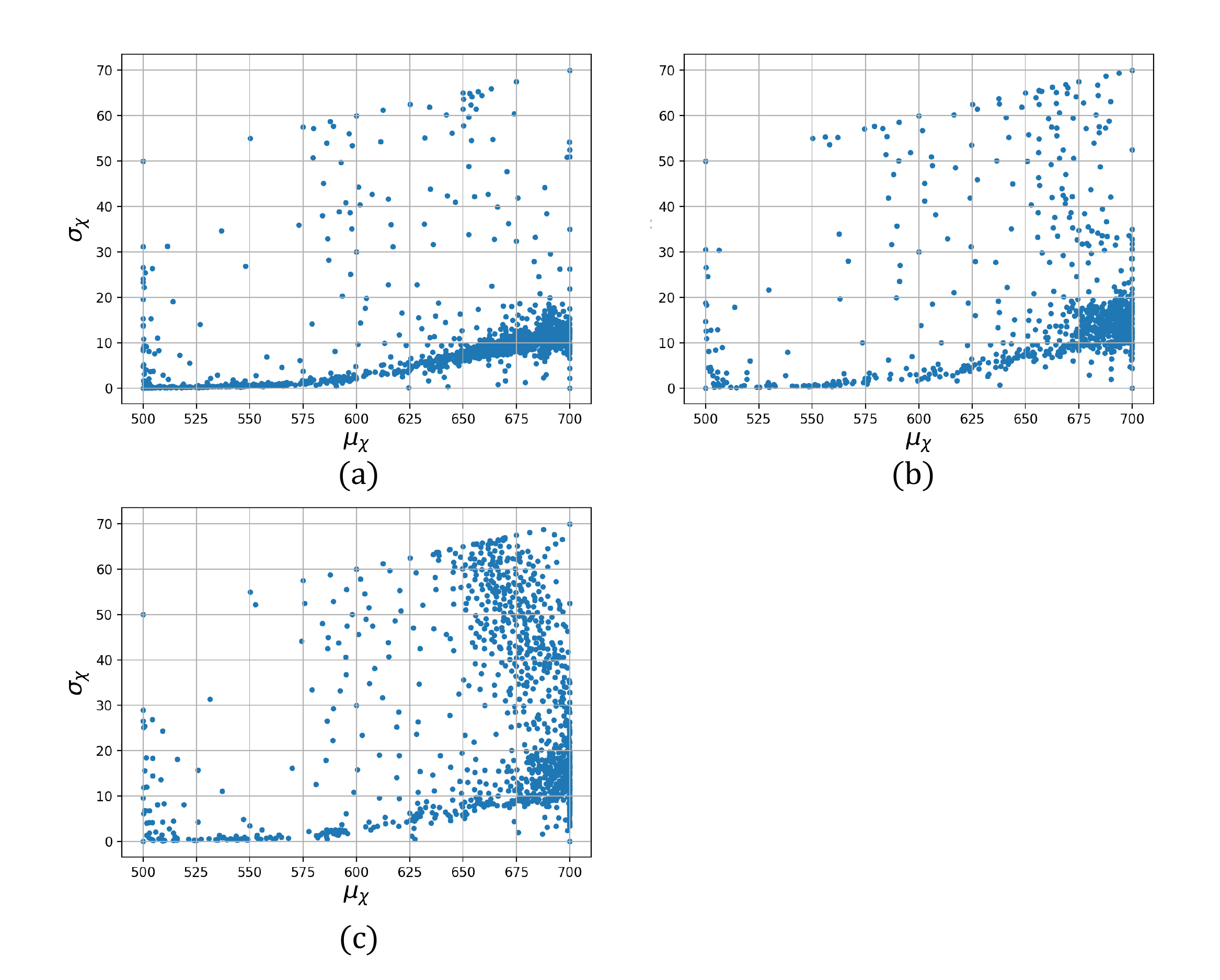}
	\caption{Final set of simulation points for different reference principal angles (a) $\theta_{\text{ref}}=\pi/10$, (b) $\theta_{\text{ref}}= \pi/12$ and (c) $\theta_{\text{ref}}=\pi/15$.}
	\label{fig:StoppingCriterion}
\end{figure}

Refinement of the DT elements continues until the convergence criterion established in Section 4 is reached. Recall that convergence is defined in terms of the average principle angle between the true solution and the solution obtained by manifold interpolation. Figure\ \ref{fig:StoppingCriterion} shows a scatter plot of the simulation points (each dot represents an STZ model simulation) at final convergence (using $d_{\mathcal{G}(\infty,\infty)}$) for three different reference principle angles ($\theta_{\text{ref}}=\frac{\pi}{10}$, $\frac{\pi}{12}$, and $\frac{\pi}{15}$). Note that the points in Figure\ \ref{fig:StoppingCriterion} are shown in the sample space of $(\mu_{\chi},\sigma_\chi)$.

This figure is illustrative of many important features of the proposed methodology and the corresponding material behavior. First, when the reference principal angle is relatively large ($\theta_{\text{ref}}=\frac{\pi}{10}$), the methodology populates the dominant transition region with a dense set of samples points (and hence a densely refined DT for interpolation). This dominant transition region corresponds to the transition from a material in which plastic deformation localizes (above the cluster of samples) and a material that deforms homogeneously (below the cluster of samples). As the reference angle is reduced ($\theta_{\text{ref}}=\frac{\pi}{10}$), a secondary transition emerges that corresponds to the transition, in a localizing material, between a single very dominant shear band (left of cluster of samples), and two less dominant shear bands (right of the cluster of samples). This is a direct consequence of the required interpolation accuracy. When the reference angle is relatively large, the interpolation on large elements is sufficiently accurate on the secondary transition but very small elements are required on the dominant transition. On the other hand, when the reference angle is small, it is necessary to highly refine both transition regions to achieve the desired level of accuracy in the interpolation. This is discussed further in the following section.

\begin{figure}[!ht]
	\centering\includegraphics[width=0.8\columnwidth]{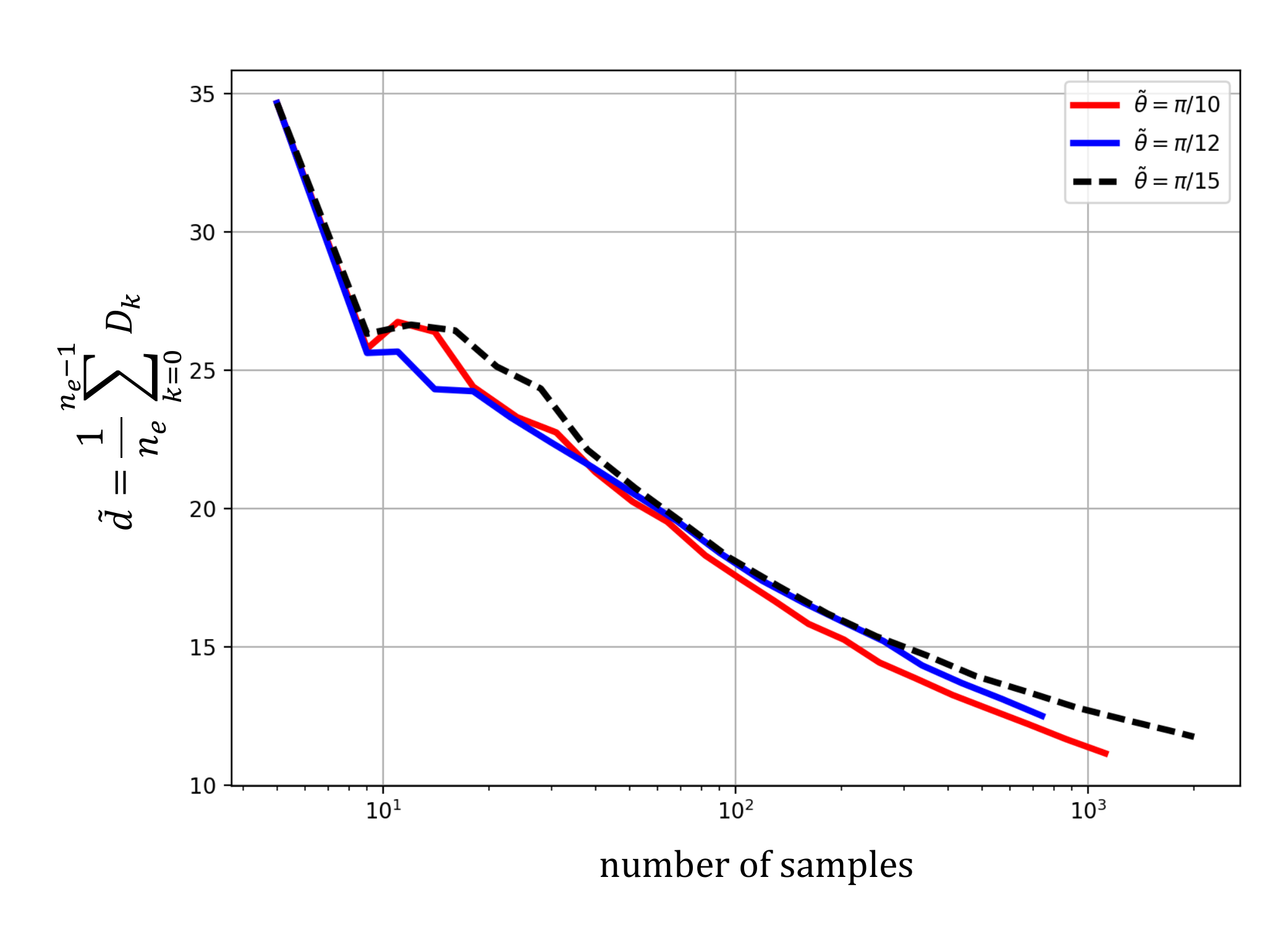}
	\caption{Convergence of Eq.\eqref{eq:34} for different reference principal angles $\theta_{\text{ref}}$ in eq. \eqref{eq:27}.}
	\label{fig:StoppingCriterion_convergence}
\end{figure}

For completeness, Figure \ref{fig:StoppingCriterion_convergence} shows that the average total Grassmann distance per element given by
\begin{equation}\label{eq:34}
\tilde{d} = \frac{1}{n_e} \sum_{k=0}^{n_e-1} D_k = \frac{1}{n_e} \sum_{k=0}^{n_e-1} \sum_{\substack{i, j=0 \\ i\neq j}}^{n_d} \delta(\boldsymbol{\Psi}_{ik},\boldsymbol{\Psi}_{jk}) 
\end{equation}
decreases exponentially with the number of samples drawn for each reference principal angle.


\subsection*{Interpolation on the manifold}

In this section, we carefully examine the performance of the interpolation. Using a set of 2000 simulations drawn according to the proposed methodology (without setting a reference principal angle for convergence criterion), we specifically study the accuracy of the interpolation in various regions of the space. 

\begin{figure}[!ht]
	\centering\includegraphics[width=1.0\columnwidth]{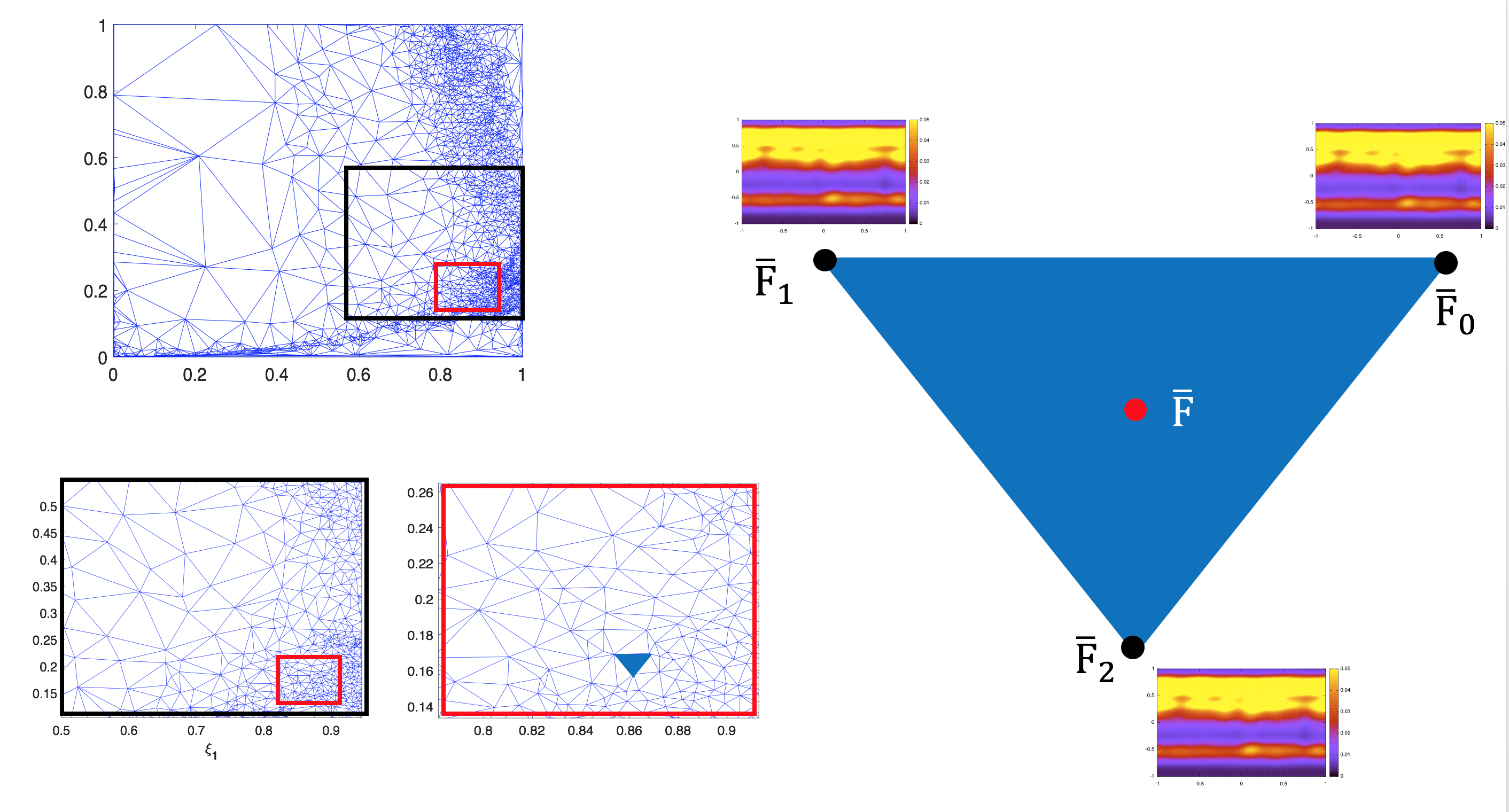}
	\caption{An element located at the bottom right area of the probability space.}
	\label{fig:area_1}
\end{figure}

\begin{figure}[!ht]
	\centering\includegraphics[width=1.0\columnwidth]{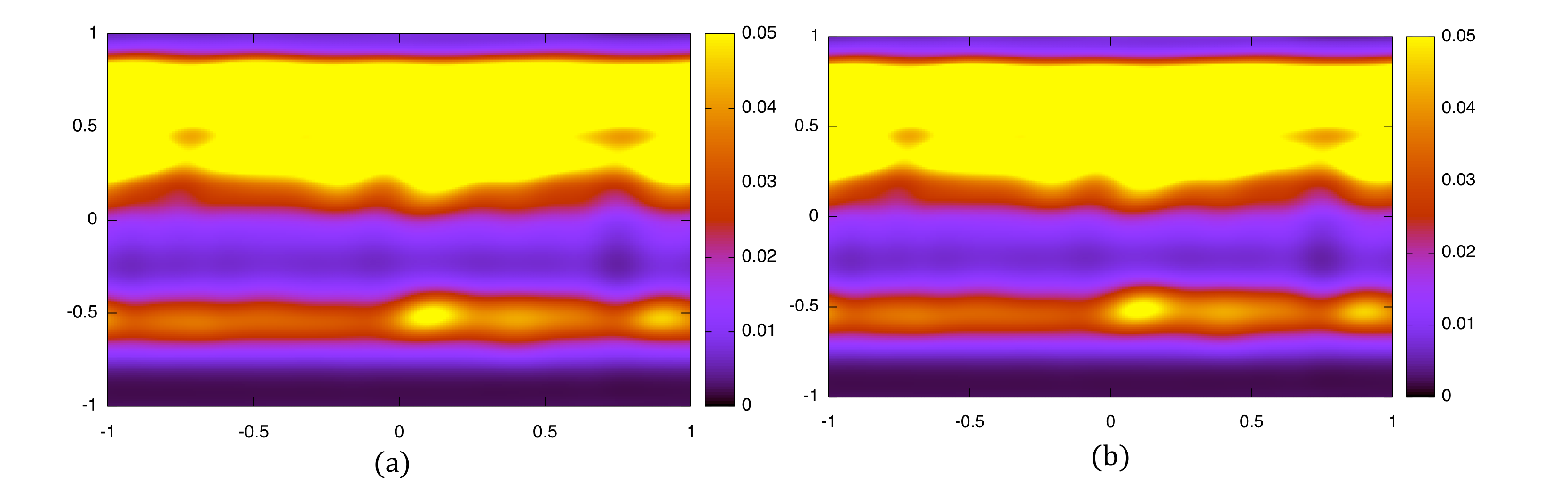}
	\caption{Plastic strain field in the interior (orange dot) of element in Fig.(\ref{fig:area_1}): (a) STZ model (b) approximated using interpolation}
	\label{fig:approximation}
\end{figure}

Fig.\ \ref{fig:area_1} gives a telescopic view of the probability space highlighting an element located in the transition region between localized and homogeneous deformation. Notice that plastic deformation exists throughout the cross-section but remains concentrated in one dominant and one secondary region. In this figure, the element (which is very small) is enlarged on the right and a point inside this element is investigated for accuracy in interpolation. 
Figure \ref{fig:approximation} shows the interpolated solution alongside the true solution from solving the STZ equations. The two are visually indistinguishable and the approximation is, indeed, quite accurate having an average principal angle of $\tilde{\theta}=0.13\approx\tfrac{\pi}{24}$.

\begin{figure}[!ht]
	\centering\includegraphics[width=1.0\columnwidth]{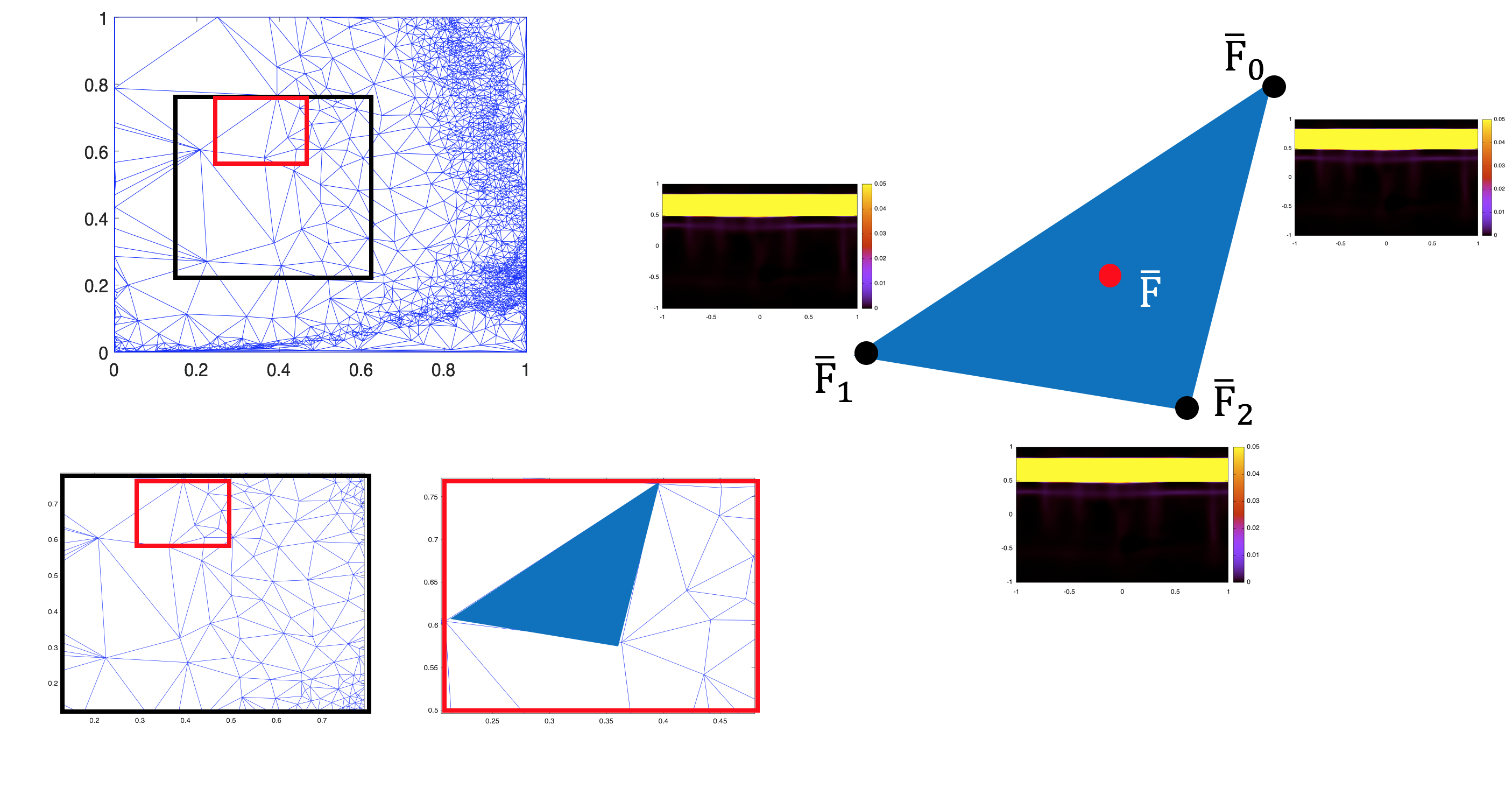}
	\caption{An element located at the top left area of the probability space.}
	\label{fig:area_2}
\end{figure}

Similarly, Figure \ref{fig:area_2} shows a telescopic view for an element located outside the transition region (in the region of strong shear band formation). Again this element (which is much larger) is enlarged on the right and a point inside is investigated for accuracy in interpolation.
Figure \ref{fig:approximation} shows the interpolated solution alongside the true solution from solving the STZ equations. The two are visually indistinguishable and the approximation is, indeed, quite accurate having an average principal angle of $\tilde{\theta}=0.13\approx\tfrac{\pi}{25}$.
\begin{figure}[!ht]
	\centering\includegraphics[width=1.0\columnwidth]{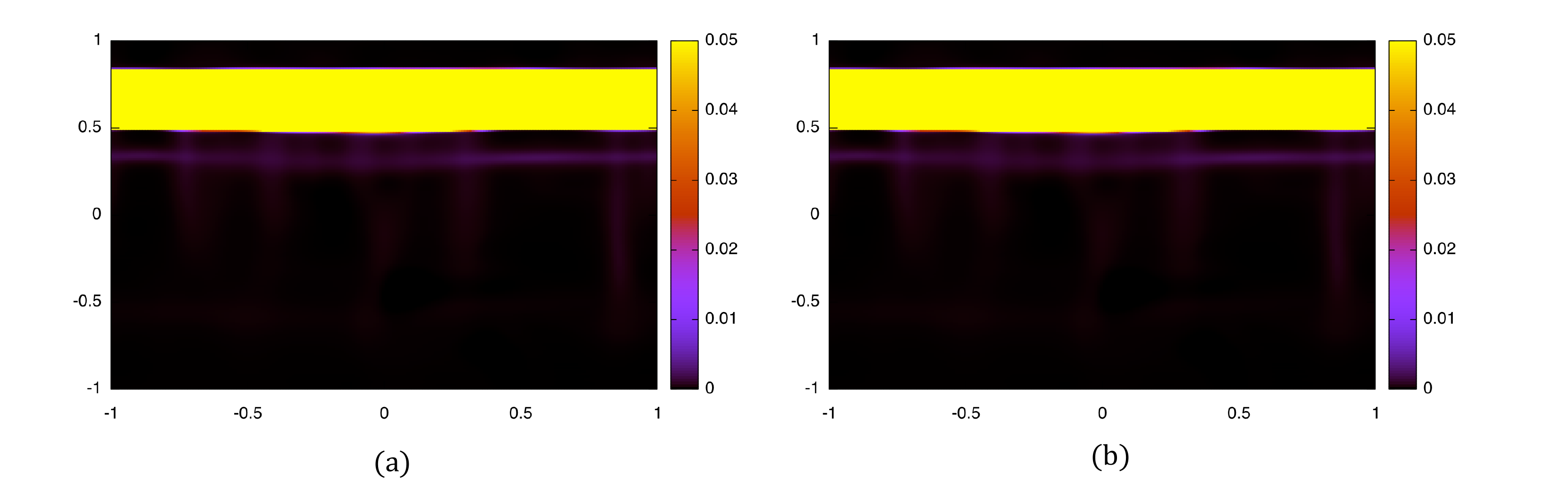}
	\caption{Plastic strain field in the interior (orange dot) of element in Fig.(\ref{fig:area_1}): (a) STZ model (b) approximated using interpolation}
	\label{fig:approximation2}
\end{figure}

Figures \ref{fig:area_1}-\ref{fig:approximation2} highlight the fact that, in regions where the behavior of the system is relatively insensitive to changes in the input parameters, interpolation will be very accurate over large regions of the parameter space (large DT elements). Meanwhile, in regions where the behavior of the systems changes rapidly with the input parameters, the interpolation requires much smaller elements in order to achieve the same level of accuracy. Note that these particular investigations were performed on a very finely resolved DT corresponding to a very high level of accuracy in the interpolation. Similar investigations for lower levels of accuracy show similar results. The desired accuracy in the interpolation remains at the users discretion and can be modified by changing the reference principal angle.

\subsubsection*{Effect of distance metric and manifold dimension}

\noindent
The proposed methodology intentionally leaves some discretion to the user in terms of assigning the manifold dimension (i.e.\ SVD truncation level) and the corresponding distance metric. In this section, we investigate these choices and show that, while a given user may have certain preferences, the methodology is robust and does not depend strongly on these choices. 

\begin{figure}[!ht]
	\centering\includegraphics[width=1.\columnwidth]{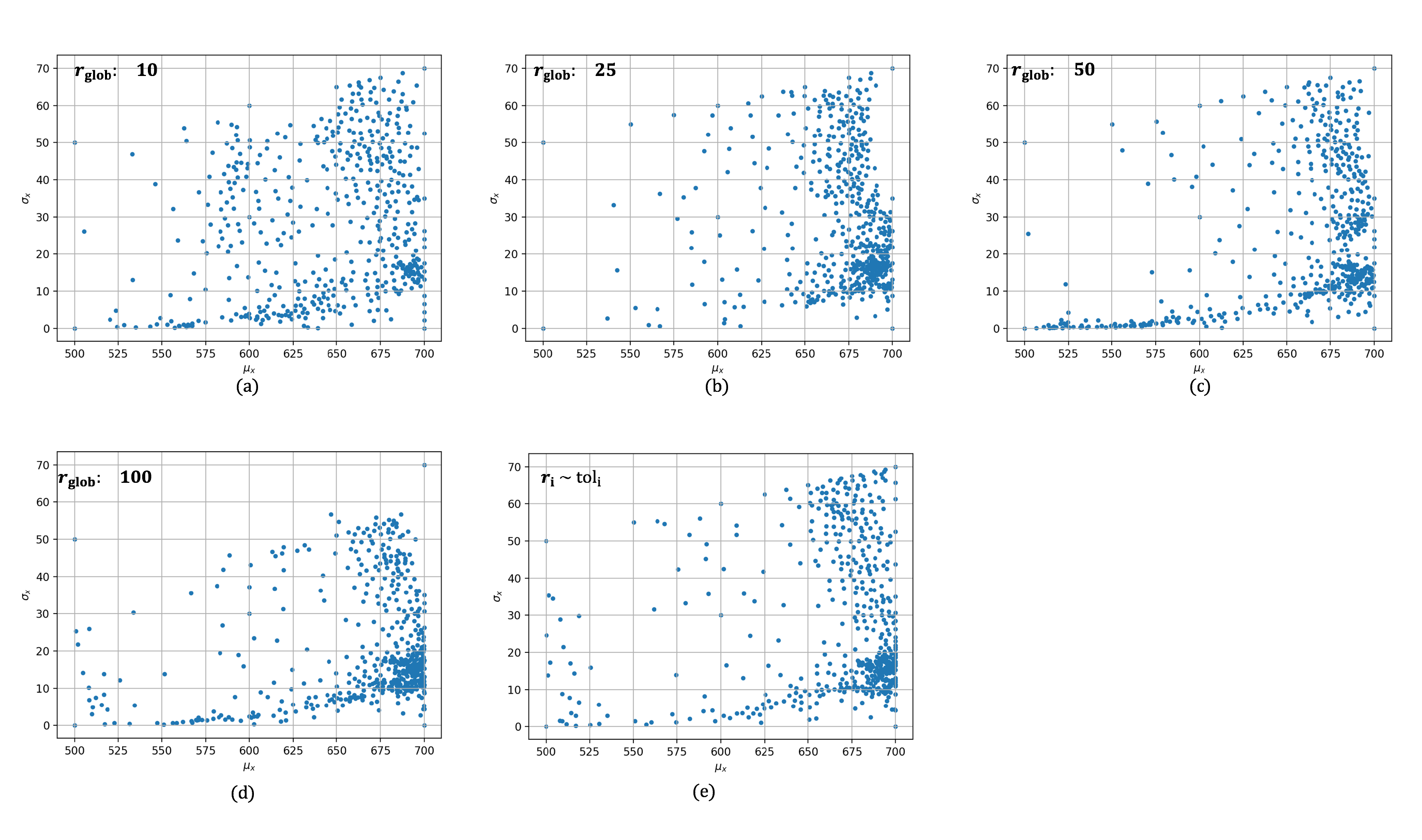}
	\caption{Samples obtained for a global rank $r_{\text{glob}}$ fixed at (a) 10, (b) 25   (c) 50 (d) 100 and (e) for an adaptive rank for each snapshot corresponding to the number of eigenvalues that are greater than a tolerance given by  $\text{tol} = \max{\Sigma}\times 201\times \text{eps}$.}
	\label{fig:samples_rank}
\end{figure}

Consider first the choice of SVD truncation level (manifold dimension). We assign a commonly used tolerance on the SVD such that each snapshot is represented through SVD with approximately equal accuracy. The consequence of this is that the dimension varies at each point requiring distances to be computed on $\mathcal{G}(\infty,\infty)$. An alternative approach is to set a fixed dimension in the SVD, $r_{\text{glob}}$, allowing distances to be computed on $\mathcal{G}(r_{\text{glob}},n)$. Figure \ref{fig:samples_rank} shows 500 sample points drawn using different values of $r_{\text{glob}}$ and those drawn using a specified tolerance level.
The sample sets show only subtle differences (most notably that the samples are more evenly distributed when only a small number of basis vectors are retained, $r_{\text{glob}}=10$) and all serve to resolve both transition regions.

\begin{figure}[!ht]
	\centering\includegraphics[width=0.7\columnwidth]{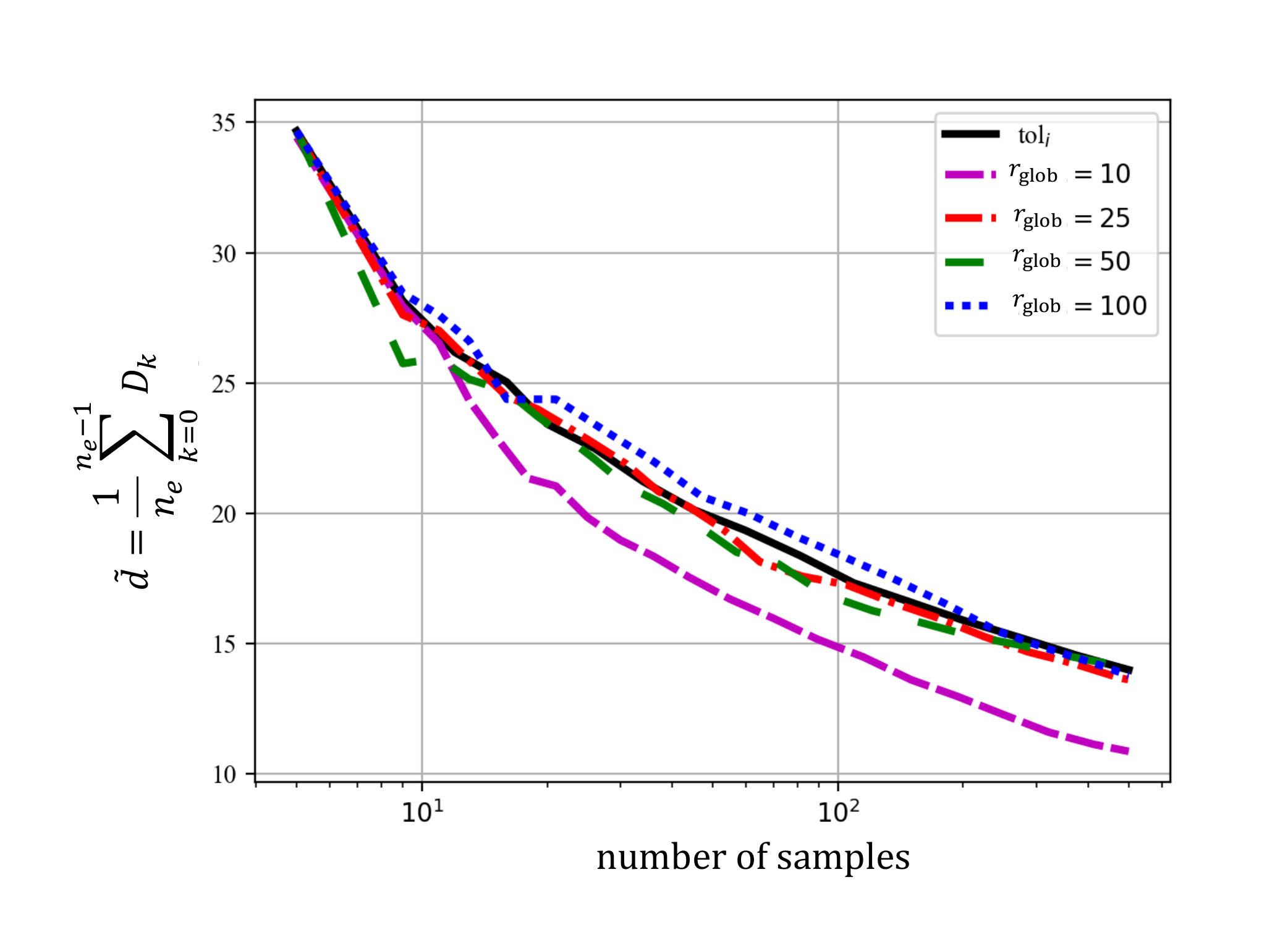}
	\caption{Convergence of the average distance $\tilde{d}$ for a global rank $r_{\text{glob}}$ fixed at values 10,  25, 50  and 100, as well as for an adaptive rank for each snapshot corresponding to the number of eigenvalues that are greater than a tolerance given by $\text{tol} = \max{\Sigma}\times n_f \times \text{eps}$.}
	\label{fig:convergence_rank}
\end{figure}

Additionally, Figure \ref{fig:convergence_rank} shows the convergence of the average elemental distance for the cases with different SVD rank (Grassmann dimension). Again, there is very little difference between the approaches with the exception again of the $r_{\text{glob}}=10$ case where distances are expected to be smaller based purely on the fact that fewer principal angles are computed.  

\begin{figure}[!ht]
	\centering\includegraphics[width=1.\columnwidth]{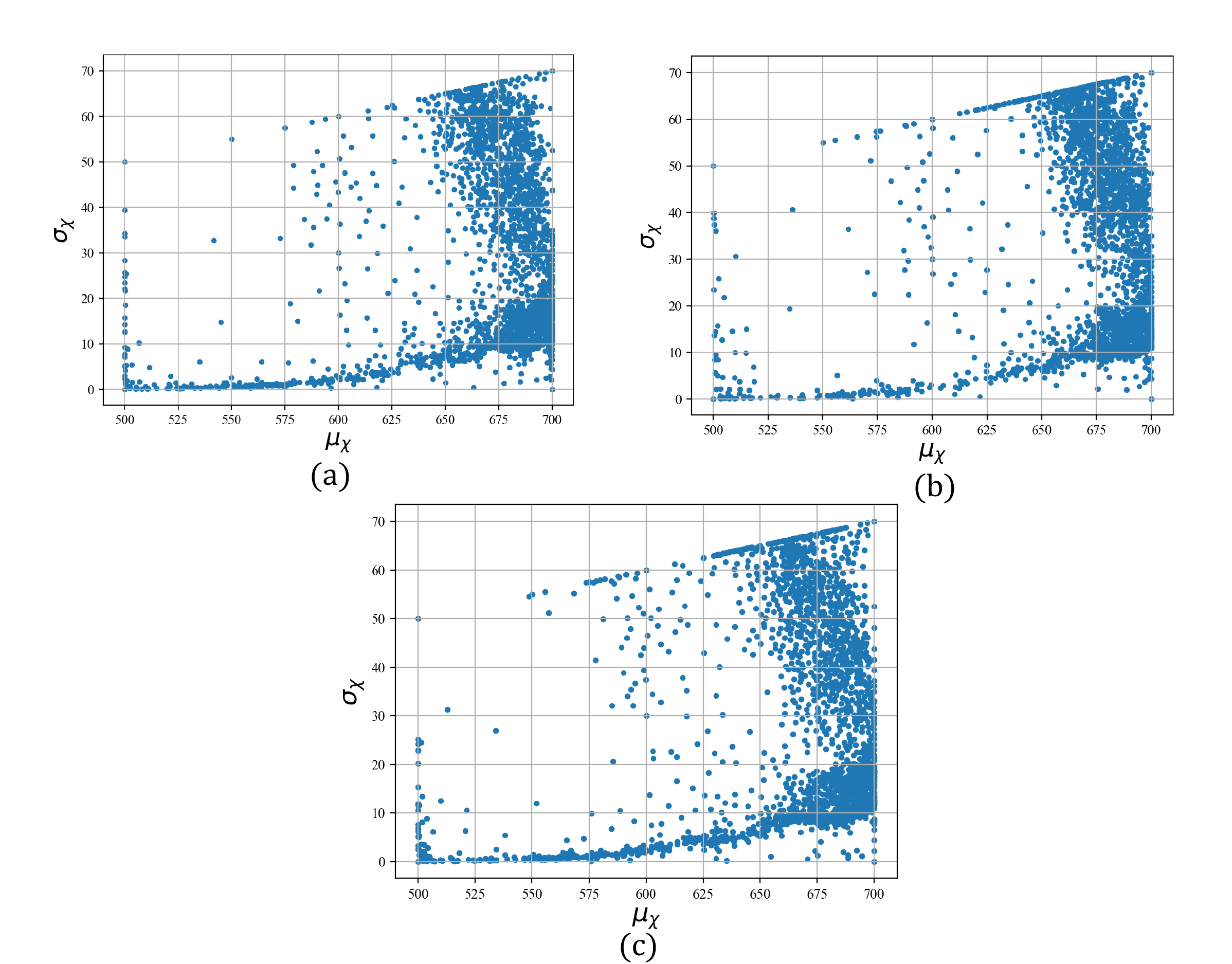}
	\caption{2000 samples obtained from the proposed method using the (a) Grassmann, (b) Chordal and (c) Procrustes metrics on $\mathcal{G}(\infty,\infty)$ (Table \ref{tab:2}).}
	\label{fig:samples}
\end{figure}

Given that the method appears insensitive to the dimension of the Grassmann manifold, we study the influence of the particular distance metric on $\mathcal{G}(\infty,\infty)$. Table \ref{tab:2} lists three different metrics (Grassmann, Chordal, and Procrustes distances) in terms of the principal angles and the relative dimensions of the manifolds at the two points. Figure \ref{fig:samples} shows 2000 samples drawn according to the proposed method using each of these different distances in the refinement criterion. 
The refinement and associated sample placement appear identical with no dependence, at all, on the distance metric. 

\begin{figure}[!ht]
	\centering\includegraphics[width=0.7\columnwidth]{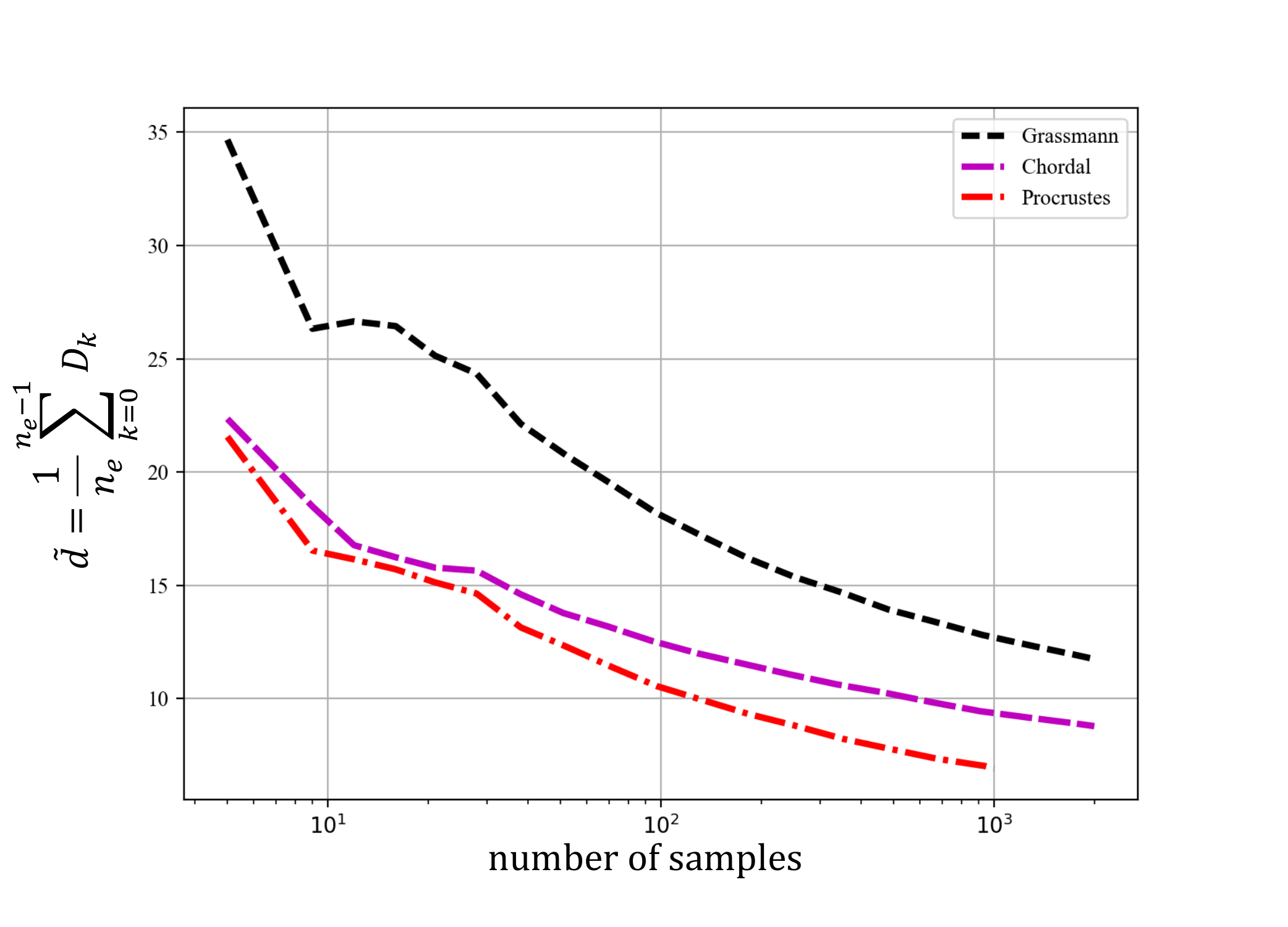}
	\caption{Average elemental distance convergence for the Grassmann, the Chordal and the Procrustes metrics on $\mathcal{G}(\infty, \infty)$.}
	\label{fig:convergence}
\end{figure}

Finally, Figure \ref{fig:convergence} shows the convergence of the average elemental distance when each of these distances is employed. Given the different scales for each metric, the magnitude of the distances are not important. The important features is that they all exhibit similar convergence trends.

\subsubsection*{Comparison with random sampling}

\begin{figure}[!ht]
	\centering
	\begin{subfigure}{0.5\columnwidth}
		\centering
		\includegraphics[scale=0.25, angle =270]{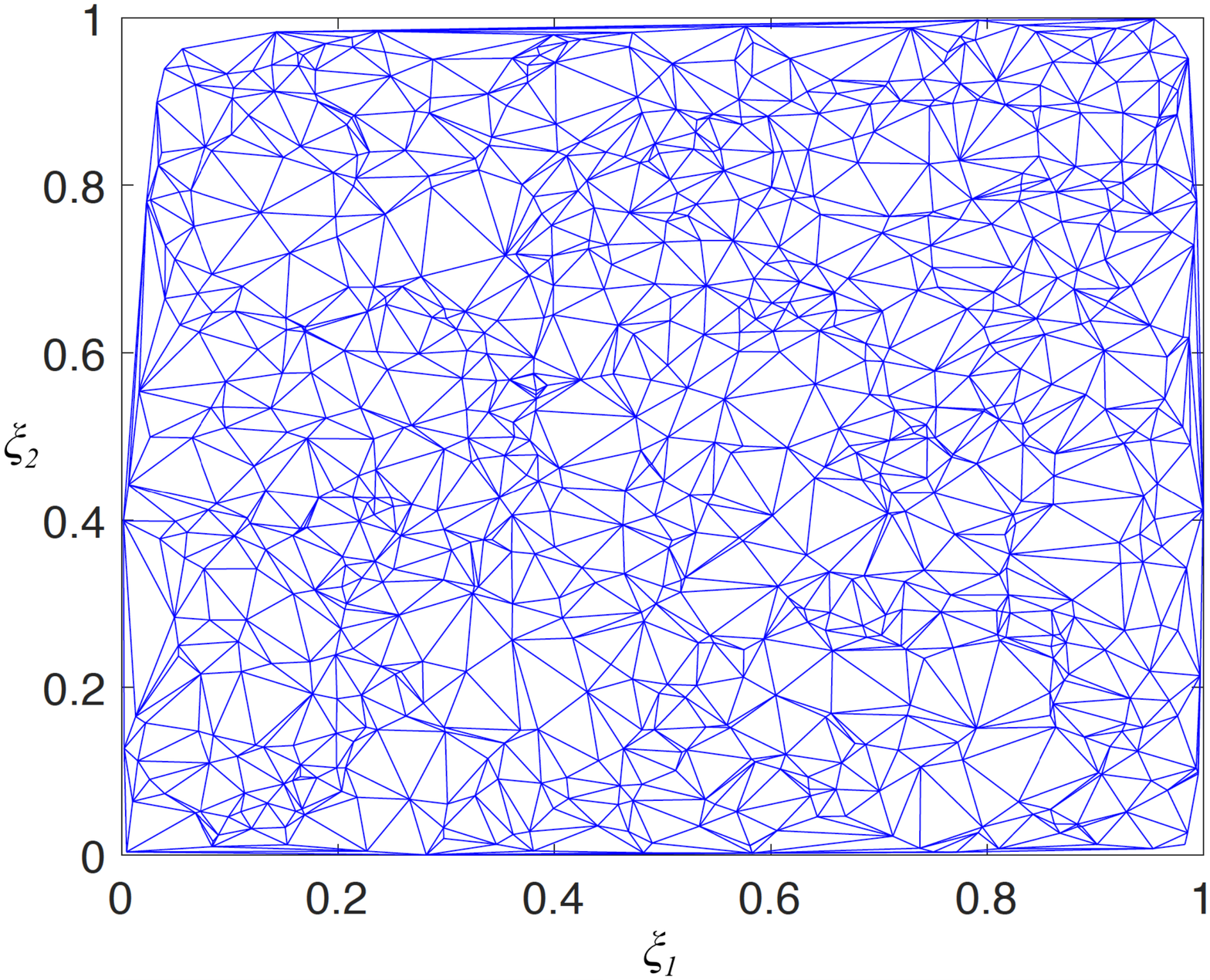}
		\caption{}
	\end{subfigure}%
	~ 
	\begin{subfigure}{0.5\columnwidth}
		\centering
		\includegraphics[scale=0.25, angle =270]{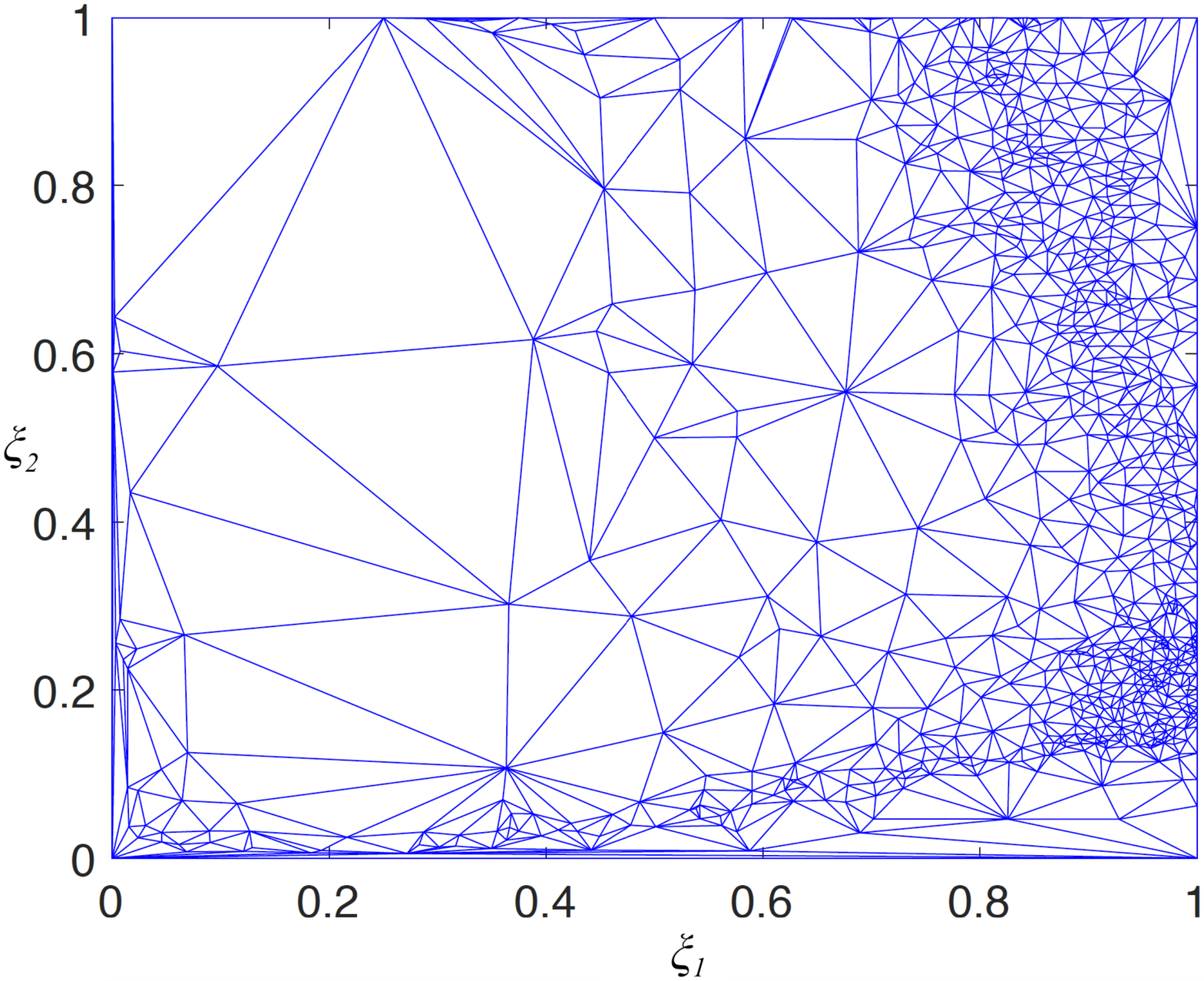}
		\caption{}
	\end{subfigure}
	\caption{Delaunay tessellations for a) random sampling and b) the proposed adaptive sampling. }
	\label{DT_compare}
\end{figure}

\noindent
In this section the accuracy of the interpolation is discussed by comparing interpolation on elements formed using the proposed adaptive sampling method and the elements formed from a set of random samples drawn uniformly in the domain.
For consistency, we utilize 800 sample points for each method and the corresponding DTs are depicted in Fig. \ref{DT_compare}

For a fair comparison of the two methods we estimate the interpolation error at the simplex centers. 
Therefore, the total number of points at which the model is approximated are 1579 and 1549 for the random sampling and the adaptive method, respectively, same as the number of simplexes in each DT. The quality of the approximations for both cases is quantified in terms of the average principal angle $\tilde{\theta}$ of Eq.\ \eqref{eq:26} and the Frobenius norm between the approximation $\tilde{\textbf{F}}$ and the true solution  $\textbf{F}$.

\begin{figure}[!ht]
\begin{subfigure}{0.5\columnwidth}
\centering
\includegraphics[scale = 0.25, angle =270]{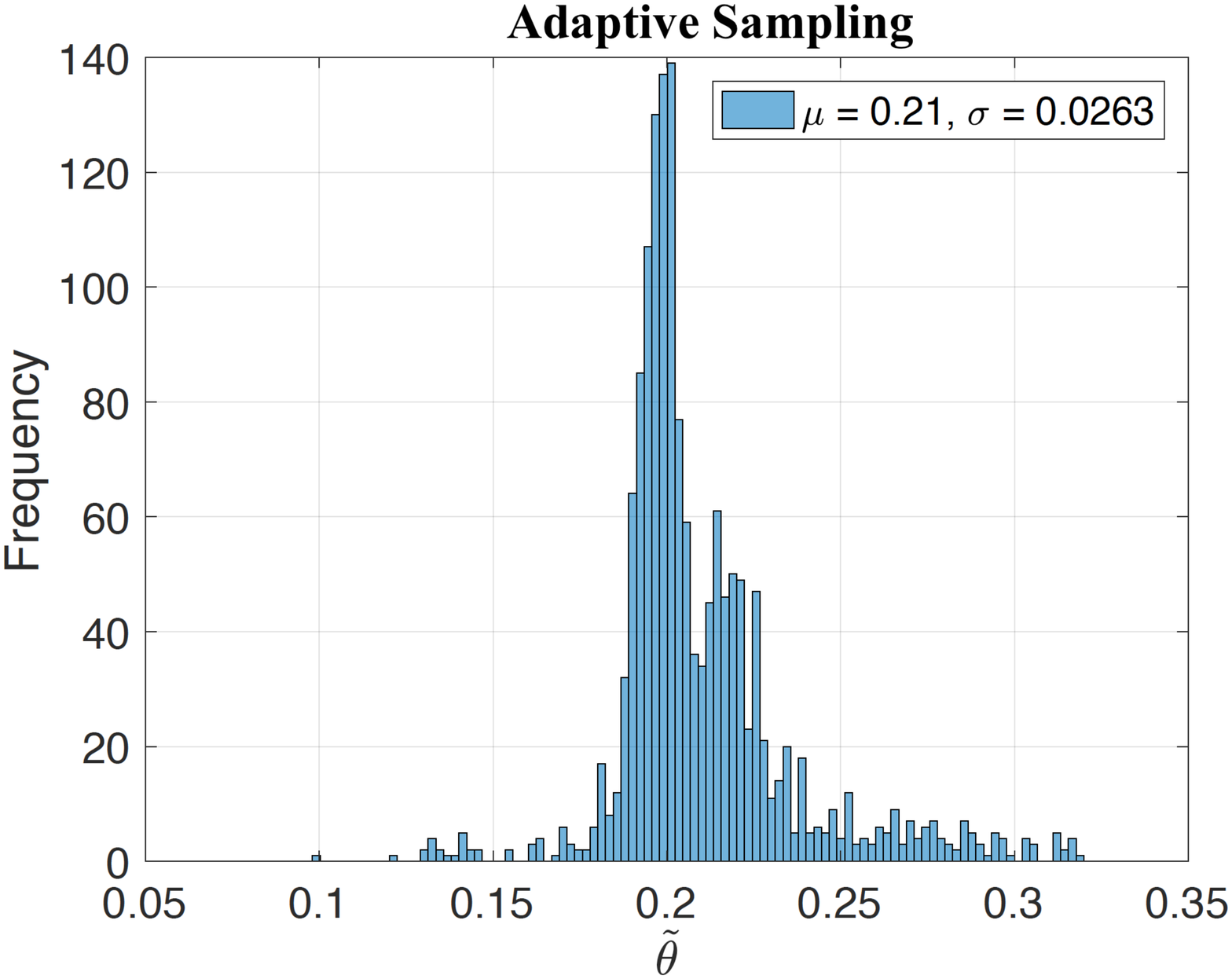}
\caption{}
\end{subfigure}
~
\begin{subfigure}{0.5\columnwidth}
\includegraphics[scale=0.25, angle =270]{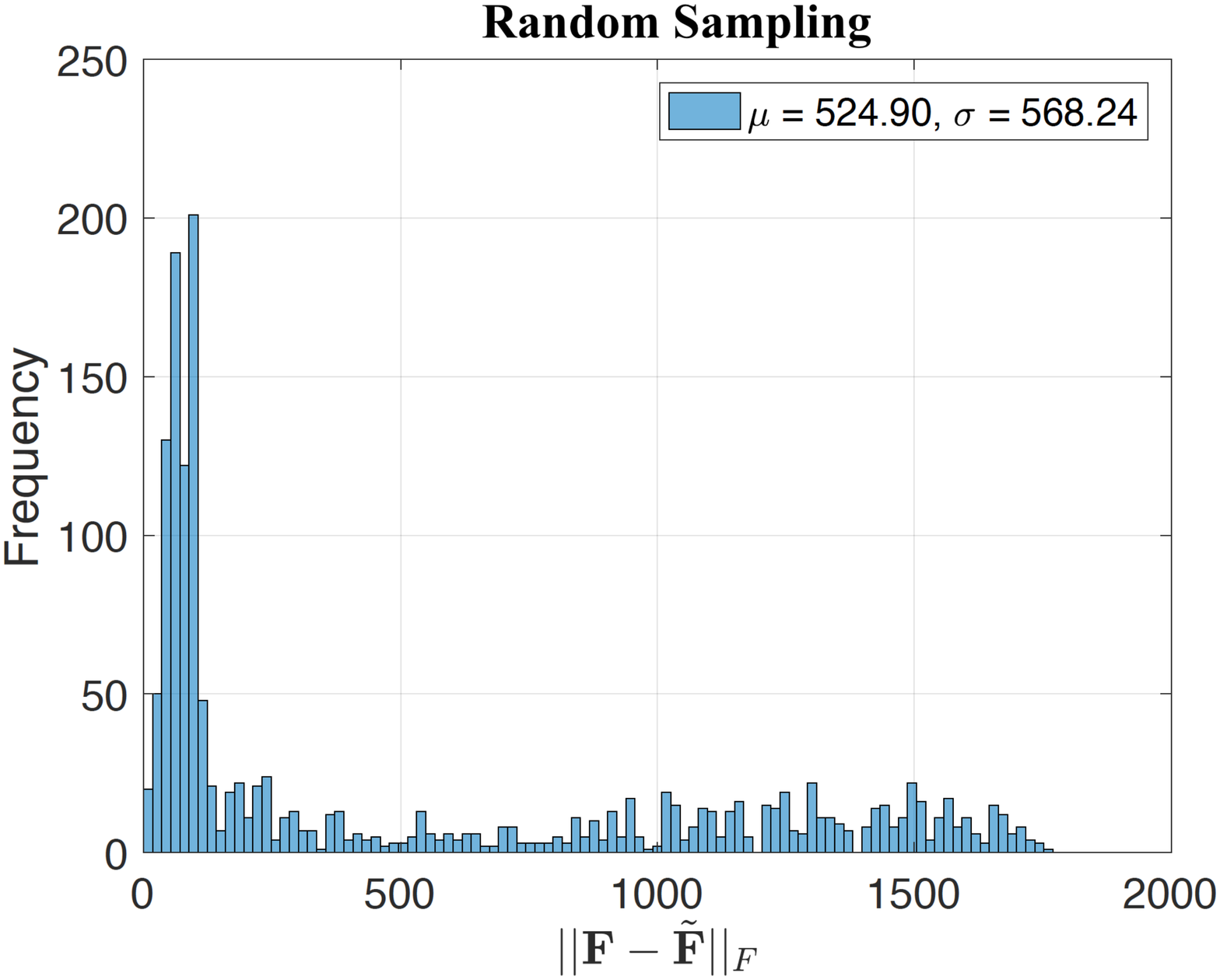}
\caption{}
\end{subfigure}

\begin{subfigure}{0.5\columnwidth}
\includegraphics[scale=0.25, angle =270]{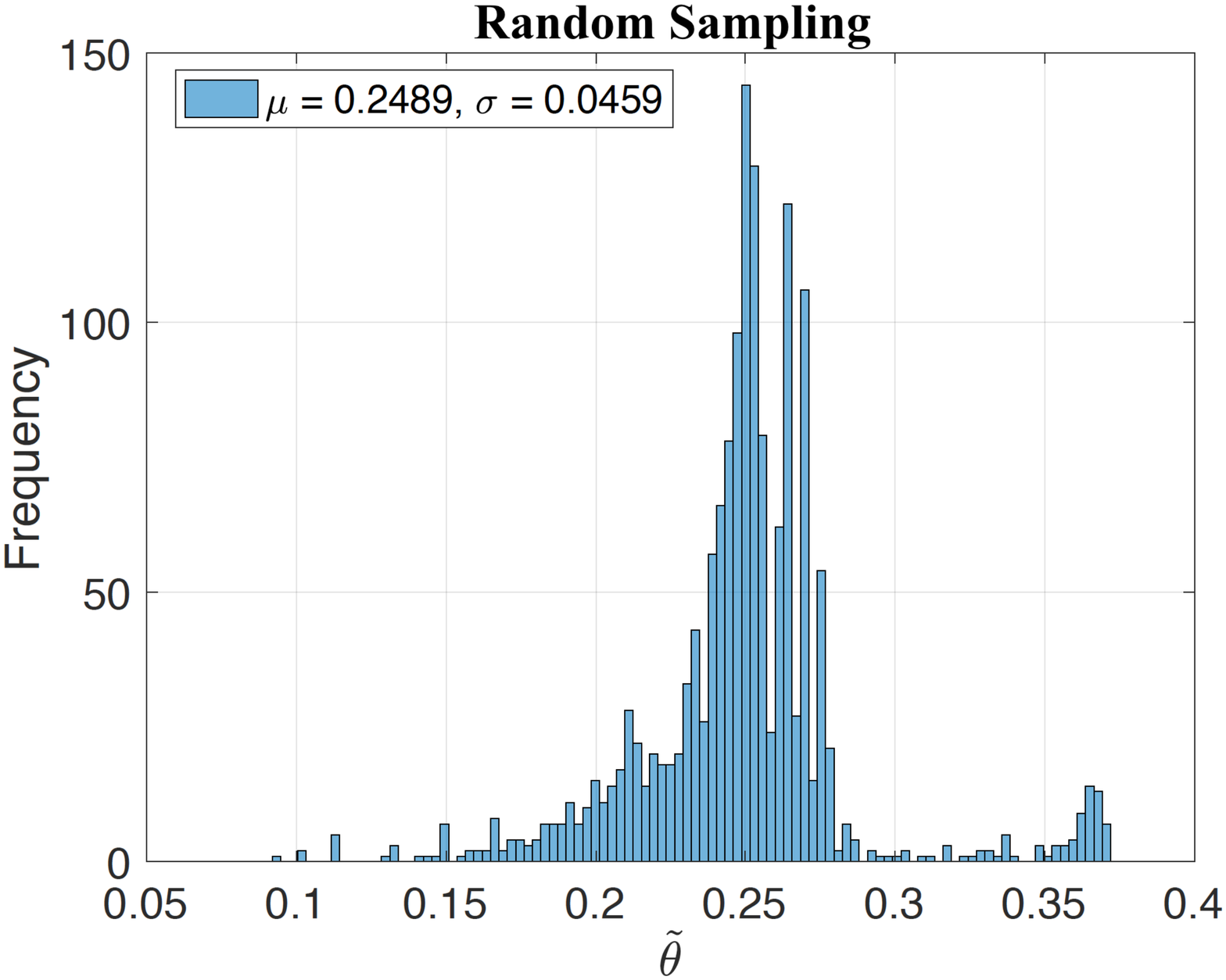}
\caption{}
\end{subfigure}
~
\begin{subfigure}{0.5\columnwidth}
\includegraphics[scale=0.25, angle =270]{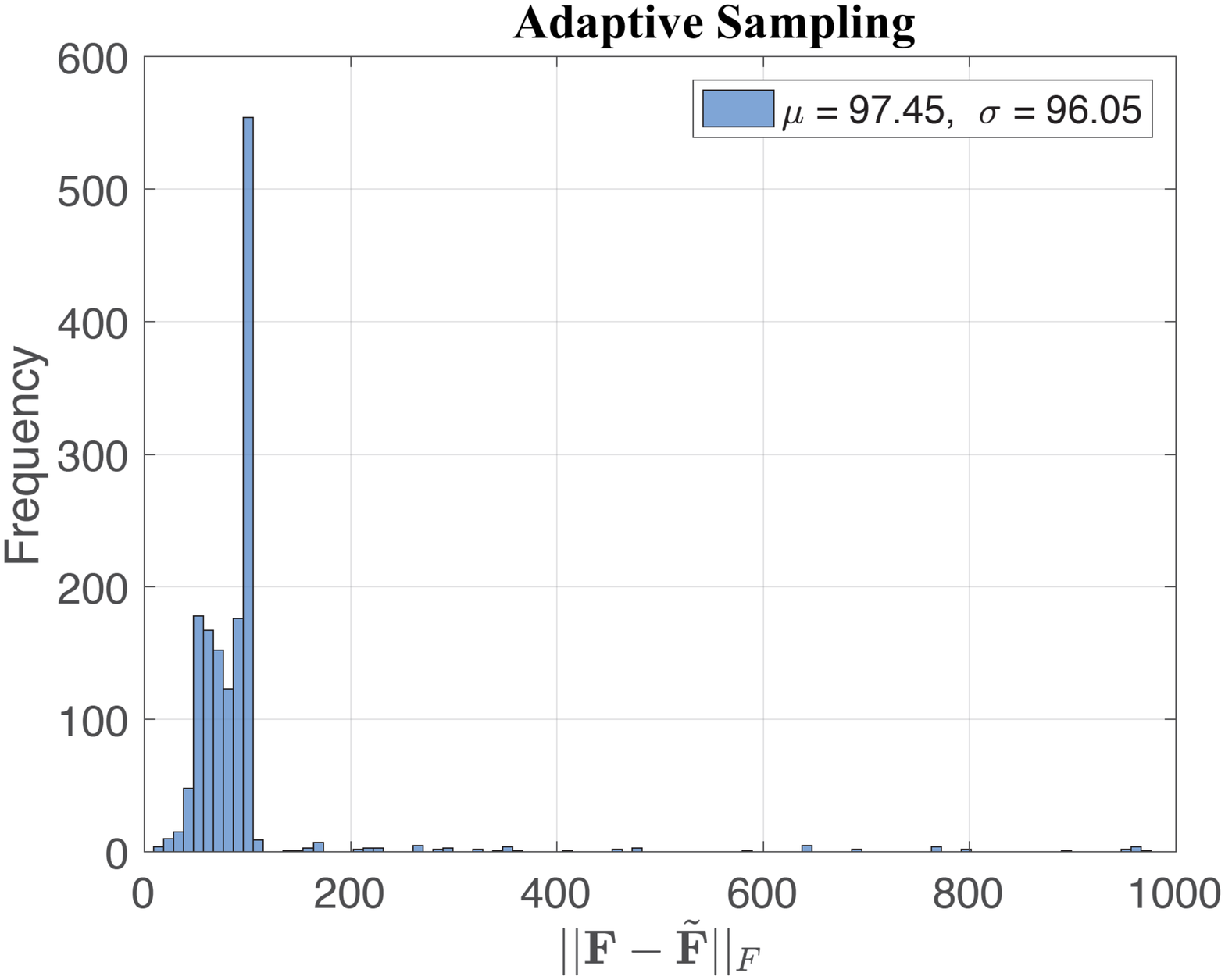}
\caption{}
\end{subfigure}
\caption{Histograms of the average principal angle and the Frobenius norm for the adaptive and non-adaptive schemes. The mean values and the standard deviations for the two metrics are also depicted in this figure. }
\label{fig:Histograms}
\end{figure}

Figure \ref{fig:Histograms} shows histograms of the average principal angle and the Frobenius norm computed at each of the simplex centers. In terms of the principal angles in Figs.\ \ref{fig:Histograms}(a) and (c), we see that interpolation from the proposed samples is consistently better than interpolation on a mesh defined by random samples. In the proposed approach the mean average principal angle is $\mu=0.21 \simeq \frac{\pi}{15}$) with standard deviation equal to $\sigma = 0.0263$ while for the random samples the mean becomes $\mu=0.25 \simeq \frac{\pi}{12.5}$) with standard deviation equal to $\sigma = 0.0459$ (nearly twice that of the proposed method). This differences in $\tilde{\theta}$ seems rather small only because the average principal angles are very sensitive. A small change in $\tilde{\theta}$ can result in a very different set of bases. With this in mind, it is worth drawing attention to the small concentration of points in Figure \ref{fig:Histograms}(c) having $\tilde{\theta}>0.35\simeq \frac{\pi}{9}$. These are very poor interpolations and the proposed method eliminates these altogether.

Perhaps a more illustrative comparison comes from the Frobenius norm, which compares the reconstructed solution instead of its Grassmann projection. These histograms are shown in Figure \ref{fig:Histograms}(b) and (d) where we see a sharp contrast between the proposed method and random sampling. The adaptive method has a strong concentration of elements where $||\mathbf{F}-\tilde{\mathbf{F}}||_F<150$ and only very few elements with large norm (all have $||\mathbf{F}-\tilde{\mathbf{F}}||_F<1000$) resulting in mean $\mu=97.45$  with standard deviation equal to $\sigma = 96.5$. The random sampling method, on the other hand, has a very large number of elements with $||\mathbf{F}-\tilde{\mathbf{F}}||_F>200$ with some even approaching $||\mathbf{F}-\tilde{\mathbf{F}}||_F=2000$ resulting in mean $\mu=524.90$ with standard deviation equal to $\sigma = 568.24$ (more than $5\times$ the proposed method). Although not shown here, our analysis indicates that those elements from the random sampling with poor interpolation lie in the regions where there is a transition in material behavior that are highly sampled by the proposed method. 

\section{Conclusions}

This paper has proposed an adaptive, stochastic simulation-based method for UQ physical systems with high-dimensional response. This is performed by projecting the solution snapshots onto a Grassmann manifold and monitoring their variations in order to decide which areas of the parameter space should be densely populated with samples. Those regions of the probability space with large distances on the Grassmann manifold correspond to regions where there are significant changes in the behavior of the system and are therefore selected for refinement. Furthermore, an approximation of the solution field is possible for a new parameter value, under the condition that the surrounding points (element vertices) are sufficiently close on the Grassman manifold. The interpolation takes place by projecting these points onto the tangent space of the manifold which is flat and thus, amenable to interpolation. The method is applied to study the probability of shear band formation in a bulk metallic glass using the shear transformation zone theory, showing the ability to identifying regions in the probability space that correspond to transitions in the material behavior, i.e.\ from homogeneous to highly localized plastic deformation.

\section*{Acknowledgements}

Methodological developments presented herein have been supported by the Office of Naval Research under grant number N-00014-16-1-2582 with Dr. Paul Hess as program officer. The work related to shear deformation of amorphous solids has been supported by the National Science Foundation Division of Materials Research under grant number 1408685 with Dr. Daryl Hess as program officer. The author is also grateful to Prof. Chris Rycroft for supplying the numerical solvers and support for the STZ equations of motion and Prof. Michael Falk for insights on the STZ plasticity theory.

\end{document}